# Neutrosophic Set is a Generalization of Intuitionistic Fuzzy Set, Inconsistent Intuitionistic Fuzzy Set (Picture Fuzzy Set, Ternary Fuzzy Set), Pythagorean Fuzzy Set (Atanassov's Intuitionistic Fuzzy Set of second type), q-Rung Orthopair Fuzzy Set, Spherical Fuzzy Set, and n-HyperSpherical Fuzzy Set, while Neutrosophication is a Generalization of Regret Theory, Grey System Theory, and Three-Ways Decision (revisited)

Florentin Smarandache
Department of Mathematics
University of New Mexico
705 Gurley Ave., Gallup, New Mexico 87301, USA

**Abstract.** In this paper we prove that Neutrosophic Set (NS) is an extension of Intuitionistic Fuzzy Set (IFS) no matter if the sum of single-valued neutrosophic components is < 1, or > 1, or = 1. For the case when the sum of components is 1 (as in IFS), after applying the neutrosophic aggregation operators one gets a different result from that of applying the intuitionistic fuzzy operators, since the intuitionistic fuzzy operators ignore the indeterminacy, while the neutrosophic aggregation operators take into consideration the indeterminacy at the same level as truth-membership and falsehood-nonmembership are taken. NS is also more flexible and effective because it handles, besides independent components, also partially independent and partially dependent components, while IFS cannot deal with these. Since there are many types of indeterminacies in our world, we can construct different approaches to various neutrosophic concepts.
Neutrosophic Set (NS) is also a generalization of Inconsistent Intuitionistic Fuzzy Set (IIFS) { which is equivalent to the Picture Fuzzy Set (PFS) and Ternary Fuzzy Set (TFS) }, Pythagorean Fuzzy Set (PyFS) {Atanassov's Intuitionistic Fuzzy Set of second type}, Spherical Fuzzy Set (SFS), n-HyperSpherical Fuzzy Set (n-HSFS), and q-Rung Orthopair Fuzzy Set (q-ROFS). And Refined Neutrosophic Set (RNS) is an extension of Neutrosophic Set. And all these sets are more general than Intuitionistic Fuzzy Set.

We prove that Atanassov's Intuitionistic Fuzzy Set of second type (AIFS2), and Spherical Fuzzy Set (SFS) do not have independent components. And we show that n-HyperSphericalFuzzy Set that we now introduce for the first time, Spherical Neutrosophic Set (SNS) and n-HyperSpherical Neutrosophic Set (n-HSNS) {the last one also introduced now for the first time} are generalizations of IFS2 and SFS.
The main distinction between Neutrosophic Set (NS) and all previous set theories are: a) the *independence* of all three neutrosophic components {truth-membership (T), indeterminacy-membership (I), falsehood-nonmembership (F)} with respect to each other in NS – while in the previous set theories their components are dependent of each other; and b) the importance of *indeterminacy* in NS - while in previous set theories indeterminacy is completely or partially ignored.

Neutrosophy is a particular case of Refined Neutrosophy, and consequently Neutrosophication is a particular case of Refined Neutrosophication. Also, Regret Theory, Grey System Theory, and Three-Ways Decision are particular cases of Neutrosophication and of Neutrosophic Probability. We have extended the Three-Ways Decision to n-Ways Decision, which is a particular case of Refined Neutrosophy.

In 2016 Smarandache defined for the first time the Refined Fuzzy Set (RFS) and Refined Fuzzy Intuitionistic Fuzzy Set (RIFS). We now, further on, define for the first time: Refined Inconsistent Intuitionistic Fuzzy Set (RIIFS){Refined Picture Fuzzy Set (RPFS), Refined Ternary Fuzzy Set (RTFS)}, Refined Pythagorean Fuzzy Set (RPyFS) {Refined Atanassov's Intuitionistic Fuzzy Set of type 2 (RAIFS2)}, Refined Spherical Fuzzy Set (RSFS), Refined n-HyperSpherical Fuzzy Set (R-n-HSFS), and Refined q-Rung Orthopair Fuzzy Set (R-q-ROFS).

**Keywords**: Neutrosophic Set, Intuitionistic Fuzzy Set, Inconsistent Intuitionistic Fuzzy Set, Picture Fuzzy Set, Ternary Fuzzy Set, Pythagorean Fuzzy Set, Atanassov's Intuitionistic Fuzzy Set of second type, Spherical Fuzzy Set, n-HyperSpherical Neutrosophic Set, q-Rung Orthopair Fuzzy Set, truth-membership, indeterminacy-membership, falsehood-nonmembership, Regret Theory, Grey System Theory, Three-Ways Decision, n-Ways Decision, Neutrosophy, Neutrosophication, Neutrosophic Probability, Refined Neutrosophy, Refined Neutrosophication.

1. Introduction

This paper recalls ideas about the distinctions between neutrosophic set and intuitionistic fuzzy set presented in previous versions of this paper [1, 2, 3, 4, 5].

Mostly, in this paper we respond to Atanassov and Vassiliev's paper [6] about the fact that neutrosophic set is a generalization of intuitionistic fuzzy set.

We use the notations employed in the neutrosophic environment [1, 2, 3, 4, 5] since they are better descriptive than the Greek letters used in intuitionistic fuzzy environment, i.e.:

truth-membership (T), indeterminacy-membership (I), and falsehood-nonmembership (F).

We also use the triplet components in this order: (T, I, F).

Neutrosophic "Fuzzy" Set (as named by Atanassov and Vassiliev [6]) is commonly called "Single-Valued" Neutrosophic Set (i.e. the neutrosophic components are single-valued numbers) by the neutrosophic community that now riches about 1,000 researchers, from 60 countries around the world, which have produced about 2,000 publications (papers, conference presentations, book chapters, books, MSc theses, and PhD dissertations).

The NS is more complex and more general than previous (crisp / fuzzy / intuitionistic fuzzy / picture fuzzy / ternary fuzzy set / Pythagorean fuzzy / Atanassov's intuitionistic fuzzy set of second type / spherical fuzzy / q-Rung orthopair fuzzy) sets, because:

- A new branch of philosophy was born, called *Neutrosophy* [7], which is a generalization of Dialectics (and of YinYang Chinese philosophy), where not only the dynamics of opposites are studied, but the dynamics of opposites together with their neutrals as well, i.e. (<A>, <neutA>, <antiA>), where <A> is an item, <antiA> its opposite, and <neutA> their neutral (indeterminacy between them).
Neutrosophy show the significance of neutrality/indeterminacy (<neutA>) that gave birth to neutrosophic set / logic / probability / statistics / measure / integral and so on, that have many practical applications in various fields.
- The sum of the Single-Valued Neutrosophic Set/Logic components was allowed to be up to 3 (this shows the importance of independence of the neutrosophic components among themselves), which permitted the characterization of paraconsistent/conflictual sets/propositions (by letting the sum of components > 1), and of paradoxical sets/propositions, represented by the neutrosophic triplet (1, 1, 1).
- NS can distinguish between *absolute* truth/indeterminacy/falsehood and *relative* truth/indeterminacy/falsehood using nonstandard analysis, which generated the Nonstandard Neutrosophic Set (NNS).
- Each neutrosophic component was allowed to take values outside of the interval [0, 1], that culminated with the introduction of the neutrosophic overset/underset/offset [8].
- NS was enlarged by Smarandache to Refined Neutrosophic Set (RNS), where each neutrosophic component was refined / split into sub-components [9]., i.e. T was refined/split into $T_1, T_2, \ldots, T_p$; I was refined / split into $I_1, I_2, \ldots, I_p$; and F was refined split into $F_1, F_2, \ldots, F_s$; where p, r, s ≥ 1 are integers and p + r + s ≥ 4; all $T_j, I_k, F_l$ are subsets of [0, 1] with no other restriction.
- RNS permitted the extension of the Law of Included Middle to the neutrosophic Law of Included Multiple-Middle [10].
- Classical Probability and Imprecise Probability were extended to Neutrosophic Probability [11], where for each event E one has: the chance that event E occurs ( ch(E) ), indeterminate-chance that event E occurs or not ( ch(neutE) ), and the chance that the event E does not occur ( ch(antiE) ), with: *0 ≤ sup{ch(E)} + sup{ch(neutE)} + sup{ch(antiE)} ≤ 3*.
- Classical Statistics was extended to Neutrosophic Statistics [12] that deals with indeterminate / incomplete / inconsistent / vague data regarding samples and populations.

And so on (see below more details). Several definitions are recalled for paper's self-containment.

## 2. Refinements of Fuzzy Types Sets

In 2016 Smarandache [8] introduced for the first time the Refined Fuzzy Set (RFS) and Refined Fuzzy Intuitionistic Fuzzy Set (RIFS).

Let $\mathcal{U}$ be a universe of discourse, and let $A \subset \mathcal{U}$ be a subset.

We give general definitions, meaning that the components may be any subsets of [0, 1]. In particular cases, the components may be single numbers, hesitant sets, intervals and so on included in [0, 1].

3. **Fuzzy Set (FS)**

$A_{FS} = \{ x(T_A(x)), x \in \mathcal{U} \}$, where $T_A: U \to \mathcal{P}([0, 1])$ is the membership degree of the generic element x with respect to the set A, and P([0, 1]) is the powerset of [0, 1], is called a Fuzzy Set.

4. **Refined Fuzzy Set (RFS)**

We have split/refined the membership degree $T_A(x)$ into sub-membership degrees. Then:

$A_{RFS} = \{x(T_A^1(x), T_A^2(x), ..., T_A^p(x)), p \geq 2, x \in U\}$, where $T_A^1(x)$ is a sub-membership degree of type 1 of the element x with respect to the set A, $T_A^2(x)$ is a sub-membership degree of type 2 of the element x with respect to the set A, ..., $T_A^p(x)$ is a sub-membership degree of type p of the element x with respect to the set A, and $T_A^j(x) \subseteq [0,1]$ for $1 \leq j \leq p$, and $\sum_{j=1}^{p} \sup T_x^j \leq 1$ for all $x \in U$.

5. **Intuitionistic Fuzzy Set (IFS)**

Let $\mathcal{U}$ be a universe of discourse, and let $A \subset \mathcal{U}$ be a subset. Then:

$A_{IFS} = \{ x(T_A(x), F_A(x)), x \in \mathcal{U} \}$, where $T_A(x), F_A(x): U \to \mathcal{P}([0, 1])$ are the membership degree respectively the nonmembership of the generic element $x$ with respect to the set $A$, and $\mathcal{P}([0, 1])$ is the powerset of [0, 1], and $\sup T_A(x) + \sup F_A(x) \leq 1$ for all $x \in U$, is called an Intuitionistic Fuzzy Set.

6. **Refined Intuitionistic Fuzzy Set (RIFS)**

We have split/refined the membership degree $T_A(x)$ into sub-membership degrees, and the nonmembership degree $F_A(x)$. Then:

$A_{RIFS} = \{x(T_A^1(x), T_A^2(x), ..., T_A^p(x); F_A^1(x), F_A^2(x), ..., F_A^s(x)), p + s \geq 3, x \in U\}$, with p, s positive nonzero integers, $\sum_{j=1}^{p} \sup T_x^j + \sum_{l=1}^{s} \sup F_x^l \leq 1$, and $T_A^j(x), F_A^l(x) \subseteq [0,1]$ for $1 \leq j \leq p$ and $1 \leq l \leq s$.

Where $T_A^1(x)$ is a sub-membership degree of type 1 of the element x with respect to the set A, $T_A^2(x)$ is a sub-membership degree of type 2 of the element x with respect to the set A, ..., $T_A^p(x)$ is a sub-membership degree of type *p* of the element x with respect to the set A.

And $F_A^1(x)$ is a sub-nonmembership degree of type 1 of the element x with respect to the set A, $F_A^2(x)$ is a sub-nonmembership degree of type 2 of the element x with respect to the set A, ..., $F_A^s(x)$ is a sub-nonmembership degree of type *s* of the element x with respect to the set A.

7. **Inconsistent Intuitionistic Fuzzy Set (IIFS) { Picture Fuzzy Set (PFS), Ternary Fuzzy Set (TFS) }**

Are defined as below:

$A_{IIFS} = A_{PFS} = A_{TFS} = \{\langle x, T_A(x), I_A(x), F_A(x) \rangle | x \in \mathcal{U}\}$,

where $T_A(x), I_A(x), F_A(x) \in P([0,1])$ and the sum $0 \leq supT_A(x) + supI_A(x) + supF_A(x) \leq 1$, for all $x \in \mathcal{U}$.

In these sets, the denominations are:

$T_A(x)$ is called *degree of membership* (or validity, or positive membership);

$I_A(x)$ is called *degree of neutral membership*;

$F_A(x)$ is called *degree of nonmembership* (or nonvalidity, or negative membership).

The *refusal degree* is: $R_A(x) = 1 - T_A(x) - I_A(x) - F_A(x) \in [0,1]$, for all $x \in \mathcal{U}$.

8. **Refined Inconsistent Intuitionistic Fuzzy Set (RIIFS) { Refined Picture Fuzzy Set (RPFS), Refined Ternary Fuzzy Set (RTFS) }**

$A_{RIIFS} = A_{RPFS} = A_{RTFS} = \{x(T_A^1(x), T_A^2(x), ..., T_A^p(x); I_A^1(x), I_A^2(x), ..., I_A^r(x); F_A^1(x), F_A^2(x), ..., F_A^s(x)), p + r + s \geq 4, x \in U\}$,

with *p, r, s* positive nonzero integers, and $T_A^j(x), I_A^k(x), F_A^l(x) \subseteq [0,1]$, for $1 \leq j \leq p$, $1 \leq k \leq r$, and $1 \leq l \leq s$, $0 \leq \sum_1^p supT_A^j(x) + \sum_1^r supI_A^k(x) + \sum_1^s supF_A^l(x) \leq 1$.

$T_A^j(x)$ is called *degree of sub-membership* (or sub-validity, or positive sub-membership) of type *j* of the element *x* with respect to the set *A*;

$I_A^k(x)$ is called *degree of sub-neutral membership* of type *k* of the element *x* with respect to the set *A*;

$F_A^l(x)$ is called *degree of sub-nonmembership* (or sub-nonvalidity, or negative sub-membership) of type *l* of the element *x* with respect to the set *A*;

and the *refusal degree* is:

$R_A(x) = [1,1] - \sum_1^p T_A^j(x) - \sum_1^r I_A^k(x) - \sum_1^s F_A^l(x) \subseteq [0,1]$, for all $x \in \mathcal{U}$.

9. **Definition of single-valued Neutrosophic Set (NS)**

Introduced by Smarandache [13, 14, 15] in 1998. Let U be a universe of discourse, and a set $A_{NS} \subseteq U$.

Then $A_{NS} = \{<x, T_A(x), I_A(x), F_A(x)> | x \in U\}$, where $T_A(x), I_A(x), F_A(x) : U \rightarrow [0, 1]$ represent the degree of truth-membership, degree of indeterminacy-membership, and degree of false-nonmembership respectively, with $0 \leq T_A(x) + I_A(x) + F_A(x) \leq 3$.

The neutrosophic components T_A(x), I_A(x), F_A(x) are <u>independent</u> with respect to each other.

### 10. Definition of single-valued Refined Neutrosophic Set (RNS)

Introduced by Smarandache [9] in 2013. Let U be a universe of discourse, and a set $A_{RNS} \subseteq$ U. Then

$A_{RNS}$ = {<x, $T_{1A}(x)$, $T_{2A}(x)$, …, $T_{pA}(x)$; $I_{1A}(x)$, $I_{2A}(x)$, …, $I_{rA}(x)$; $F_{1A}(x)$, $F_{2A}(x)$, …, $F_{sA}(x)$> | x ∈ U}, where all $T_{jA}(x)$, $1 \le j \le p$, $I_{kA}(x)$, $1 \le k \le r$, $F_{lA}(x)$, $1 \le l \le s$, : U → [0, 1], and

$T_{jA}(x)$ represents the j-th sub-membership degree,

$I_{kA}(x)$ represents the k-th sub-indeterminacy degree,

$F_{lA}(x)$ represents the *l*-th sub-nonmembership degree,

with $p, r, s \ge 1$ integers, where $p + r + s = n \ge 4$, and:

$$0 \le \sum_{j=1}^{p} T_{jA}(x) + \sum_{k=1}^{r} I_{kA}(x) + \sum_{l=1}^{s} T_{jA}(x) \le n.$$

All neutrosophic sub-components $T_{jA}(x)$, $I_{kA}(x)$, $F_{lA}(x)$ are <u>independent</u> with respect to each other.

Refined Neutrosophic Set is a generalization of Neutrosophic Set.

### 11. Definition of single-valued Intuitionistic Fuzzy Set (IFS)

Introduced by Atanassov [16, 17, 18] in 1983. Let U be a universe of discourse, and a set $A_{IFS} \subseteq$ U. Then $A_{IFS}$ = {<x, $T_A(x)$, $F_A(x)$> | x ∈ U}, where $T_A(x)$, $F_A(x)$ : U → [0, 1] represent the degree of membership and degree of nonmembership respectively, with $T_A(x) + F_A(x) \le 1$, and $I_A(x) = 1 - T_A(x) - F_A(x)$ represents degree of indeterminacy (in previous publications it was called degree of hesitancy).

The intuitioinistic fuzzy components $T_A(x)$, $I_A(x)$, $F_A(x)$ are <u>dependent</u> with respect to each other.

### 12. Definition of single-valued Inconsistent Intuitionistic Fuzzy Set (equivalent to single-valued Picture Fuzzy Set, and with single-valued Ternary Fuzzy Set)

The single-valued *Inconsistent Intuitionistic Fuzzy Set* (IIFS), introduced by Hindde and Patching [19] in 2008, and the single-valued *Picture Fuzzy Set* (PFS), introduced by Cuong [20] in 2013, indeed coincide, as Atanassov and Vassiliev have observed; also we add that single-valued *Ternary Fuzzy Set*, introduced by Wang, Ha and Liu [21] in 2015 also coincide with them. All these three notions are defined as follows.

Let $\mathcal{U}$ be a universe of discourse, and let's consider a subset $A \subseteq \mathcal{U}$.

Then $A_{IIFS} = A_{PFS} = A_{TFS} = \{\langle x, T_A(x), I_A(x), F_A(x) \rangle | x \in \mathcal{U}\}$,

where $T_A(x), I_A(x), F_A(x) \in [0, 1]$, and the sum $0 \le T_A(x) + I_A(x) + F_A(x) \le 1$, for all $x \in \mathcal{U}$.

In these sets, the denominations are:

$T_A(x)$ is called *degree of membership* (or validity, or positive membership);

$I_A(x)$ is called *degree of neutral membership*;

$F_A(x)$ is called *degree of nonmembership* (or nonvalidity, or negative membership).

The *refusal degree* is: $R_A(x) = 1 - T_A(x) - I_A(x) - F_A(x) \in [0, 1]$, for all $x \in \mathcal{U}$.

The IIFS (PFS, TFS) components $T_A(x), I_A(x), F_A(x), R_A(x)$ are <u>dependent</u> with respect to each other.

Wang, Ha and Liu's [21] assertion that "neutrosophic set theory is difficult to handle the voting problem, as the sum of the three components is greater than 1" is not true, since the sum of the three neutrosophic components is not necessarily greater than 1, but it can be less than or equal to any number between 0 and 3, i.e. $0 \leq T_A(x) + I_A(x) + F_A(x) \leq 3$, so for example the sum of the three neutrosophic components can be less than 1, or equal to 1, or greater than 1 depending on each application.

### 13. Inconsistent Intuitionistic Fuzzy Set and the Picture Fuzzy Set and Ternary Fuzzy Set are particular cases of the Neutrosophic Set

The *Inconsistent Intuitionistic Fuzzy Set* and the *Picture Fuzzy Set* and *Ternary Fuzzy Set* are particular cases of the *Neutrosophic Set* (NS). Because, in neutrosophic set, similarly taking single-valued components $T_A(x), I_A(x), F_A(x) \in [0, 1]$, one has the sum $T_A(x) + I_A(x) + F_A(x) \leq 3$, which means that $T_A(x) + I_A(x) + F_A(x)$ can be equal to or less than any number between 0 and 3.

Therefore, in the particular case when choosing the sum equal to $1 \in [0, 3]$ and getting $T_A(x) + I_A(x) + F_A(x) \leq 1$, one obtains IIFS and PFS and TFS.

### 14. Single-valued Intuitionistic Fuzzy Set is a particular case of single-valued Neutrosophic Set

Single-valued Intuitionistic Fuzzy Set is a particular case of single-valued Neutrosophic Set, because we can simply choose the sum to be equal to 1:

$T_A(x) + I_A(x) + F_A(x) = 1.$

### 15. Inconsistent Intuitionistic Fuzzy Set and Picture Fuzzy Set and Ternary Fuzzy Set are also particular cases of single-valued Refined Neutrosophic Set.

The *Inconsistent Intuitionistic Fuzzy Set* (IIFS), *Picture Fuzzy Set* (PFS), and *Ternary Fuzzy Set* (TFS), that coincide with each other, are in addition particular case(s) of *Single-Valued Refined Neutrosophic Set* (RNS).

We may define:

$A_{IIFS} \equiv A_{PFS} = A_{TFS} = \{x, T_A(x), I_{1_A}(x), I_{2_A}(x), F_A(x) | x \in \mathcal{U}\},$

with $T_A(x), I_{1_A}(x), I_{2_A}(x), F_A(x) \in [0, 1]$,

and the sum $T_A(x) + I_{1_A}(x) + I_{2_A}(x) + F_A(x) = 1$, for all $x \in \mathcal{U}$;

where:

$T_A(x)$ is the degree of positive membership (validity, etc.);

$I_{1_A}$ is the degree of neutral membership;

$I_{2_A}(x)$ is the refusal degree;

$F_A(x)$ is the degree of negative membership (non-validity, etc.).

$n = 4$, and as a particular case of the sum $T_A(x) + I_{1_A}(x) + I_{2_A}(x) + F_A(x) \leq 4$, where the sum can be any positive number up to 4, we take the positive number 1 for the sum:

$T_A(x) + I_{1_A}(x) + I_{2_A}(x) + F_A(x) = 1.$

## 16. Independence of Neutrosophic Components vs. Dependence of Intuitionistic Fuzzy Components

Section 4, equations (46) - (51) in Atanassov's and Vassiliev's paper [6] is reproduced below:

> *"4. Interval valued intuitionistic fuzzy sets, intuitionistic fuzzy sets, and neutrosophic fuzzy sets*
> 
> (…) the concept of a Neutrosophic Fuzzy Set (NFS) is introduced, as follows:
> 
> $$A^n = \{x, \mu_A^n(x), v_A^n(x), \pi_A^n(x) | x \in E\}, \quad (46)$$
> 
> where $\mu_A^n(x), v_A^n(x), \pi_A^n(x) \in [0,1]$, and have the same sense as IFS.
> 
> Let
> 
> $$\sup_{y \in E} \mu_A^n(y) + \sup_{y \in E} v_A^n(y) + \sup_{y \in E} \pi_A^n(y) \neq 0. \quad (47)$$
> 
> Then we define:
> 
> $$\mu_A^i(x) = \frac{\mu_A^n(x)}{\sup_{y \in E} \mu_A^n(y) + \sup_{y \in E} v_A^n(y) + \sup_{y \in E} \pi_A^n(y)}; \quad (48)$$
> 
> $$v_A^i(x) = \frac{v_A^n(x)}{\sup_{y \in E} \mu_A^n(y) + \sup_{y \in E} v_A^n(y) + \sup_{y \in E} \pi_A^n(y)}; \quad (49)$$
> 
> $$\pi_A^i(x) = \frac{\pi_A^n(x)}{\sup_{y \in E} \mu_A^n(y) + \sup_{y \in E} v_A^n(y) + \sup_{y \in E} \pi_A^n(y)}; \quad (50)$$
> 
> $$i_A^i(x) = 1 - \mu_A^i(y) - v_A^i(y) - \pi_A^i(y). \quad (51)"$$

Using the neutrosophic component common notations, $T_A(x) \equiv \mu_A^n(x)$, $I_A(x) \equiv \pi_A^n(x)$, and $F_A(x) \equiv v_A^n(x)$, the refusal degree $R_A(x)$, and $A_N \equiv A^n$ for the neutrosophic set, and considering the triplet's order (T, I, F), with the universe of discourse $\mathcal{U} \equiv E$, we can re-write the above formulas as follows:

$$A_N = \{\langle x_1, T_A(x), I_A(x), T_A(x)\rangle | x \in \mathcal{U}\}, \tag{46)'}$$

where $T_A(x), I_A(x), F_A(x) \in [0, 1]$, for all $x \in \mathcal{U}$.

Neutrosophic *Fuzzy* Set is commonly named *Single-Valued* Neutrosophic Set (SVNS), i.e. the components are single-valued numbers.

The authors, Atanassov and Vassiliev, assert that $T_A(x), I_A(x), F_A(x)$ "have the same sense as IFS" (Intuitionistic Fuzzy Set).

But this is untrue, since in IFS one has $T_A(x) + I_A(x) + F_A(x) \leq 1$, therefore the IFS components $T_A(x), I_A(x), T_A(x)$ are <u>dependent</u>, while in SVNS (Single-Valued Neutrosophic Set), one has $T_A(x) + I_A(x) + F_A(x) \leq 3$, what the authors omit to mention, therefore the SVNS components $T_A(x), I_A(x), F_A(x)$ are <u>independent</u>, and this makes a big difference, as we'll see below.

In general, for the *dependent components*, if one component's value changes, the other components values also change (in order for their total sum to keep being up to 1). While for the *independent components*, if one component changes, the other components do not need to change since their total sum is always up to 3.

Let's re-write the equations (47) - (51) from authors' paper:

Assume

$$\sup_{y\in\mathcal{U}} T_A(y) + \sup_{y\in\mathcal{U}} I_A(y) + \sup_{y\in\mathcal{U}} F_A(x) \neq 0. \tag{47)'}$$

The authors have defined:

$$T_A^{IIFS}(x) = \frac{T_A(x)}{\sup_{y\in\mathcal{U}} T_A(y) + \sup_{y\in\mathcal{U}} I_A(y) + \sup_{y\in\mathcal{U}} F_A(y)}; \tag{48)'}$$

$$I_A^{IIFS}(x) = \frac{I_A(x)}{\sup_{y\in\mathcal{U}} T_A(y) + \sup_{y\in\mathcal{U}} I_A(y) + \sup_{y\in\mathcal{U}} F_A(y)}; \tag{50)'}$$

$$F_A^{IIFS}(x) = \frac{F_A(x)}{\sup_{y\in\mathcal{U}} T_A(y) + \sup_{y\in\mathcal{U}} I_A(y) + \sup_{y\in\mathcal{U}} F_A(y)}. \tag{49)'}$$

These mathematical transfigurations, which transform [change in form] the neutrosophic components $T_A(x), I_A(x), F_A(x) \in [0, 1]$, whose sum

$T_A(x) + I_A(x) + F_A(x) \leq 3$, into inconsistent intuitionistic fuzzy components:

$T_A^{IIFS}(x), I_A^{IIFS}(x), F_A^{IIFS}(x) \in [0, 1]$,

whose sum $T_A^{IIFS}(x) + I_A^{IIFS}(x) + F_A^{IIFS}(x) \leq 1$,

and the refusal degree

$$R_A^{IIFS}(x) = 1 - T_A^{IIFS}(x) - I_A^{IIFS}(x) - F_A^{IIFS}(x) \in [0,1], \qquad (51)'$$

distort the original application, i.e. the original neutrosophic application and its intuitioinistic fuzzy transformed application are not equivalent, see below.

This is because, in this case, *the change in form brings a change in content*.

### 17. By Transforming the Neutrosophic Components into Intuitionistic Fuzzy Components the Independence of the Neutrosophic Components is Lost

In reference paper [6], Section 4, Atanassov and Vassilev convert the neutrosophic components into intuitionistic fuzzy components.
But, converting a single-valued neutrosophic triplet $(T_1, I_1, F_1)$, with $T_1, I_1, F_1 \in [0, 1]$ and $T_1 + I_1 + F_1 \leq 3$ that occurs into a neutrosophic application $\alpha_N$, to a single-valued intuitionistic triplet $(T_2, I_2, F_2)$, with $T_2, I_2, F_2 \in [0, 1]$ and $T_2 + F_2 \leq 1$ (or $T_2 + I_2 + F_2 = 1$) that would occur into an intuitionistic fuzzy application $\alpha_{IF}$, is just a *mathematical artifact*, and there could be constructed many such mathematical operators [the authors present four of them], even more: it is possible to convert from the sum $T_1 + I_1 + F_1 \leq 3$ to the sum
$T_2 + I_2 + F_2$ equals to any positive number – but they are just *abstract transformations*.
The neutrosophic application $\alpha_N$ will not be equivalent to the resulting intuitionistic fuzzy application $\alpha_{IF}$, since while in $\alpha_N$ the neutrosophic components $T_1, I_1, F_1$ are <u>independent</u> (because their sum is up to 3), in $\alpha_{IF}$ the intuitionistic fuzzy components $T_2, I_2, F_2$ are <u>dependent</u> (because their sum is 1). Therefore, the independence of components is lost.
And the *independence* of the neutrosophic components is the main distinction between neutrosophic set vs. intuitionistic fuzzy set.
Therefore, the resulted intuitionistic fuzzy application $\alpha_{IF}$ after the mathematical transformation is just a subapplication (particular case) of the original neutrosophic application $\alpha_N$.

### 18. Degree of Dependence/Independence between the Components

The degree of dependence/independence between components was introduced by Smarandache [22] in 2006.

In general, the sum of two components x and y that vary in the unitary interval [0, 1] is:

*0 ≤ x+y ≤ 2-d(x,y)*, where *d(x,y)* is the degree of dependence between *x* and *y*, while *1-d(x,y)* is the degree of independence between x and y.

NS is also flexible because it handles, besides independent components, also partially independent and partially dependent components, while IFS cannot deal with these.

For example, if *T* and *F* are totally dependent, then *0 ≤ T + F ≤ 1*, while if component *I* is independent from them, thus *0 ≤ I ≤ 1*, then *0 ≤ T + I + F ≤ 2*. Therefore the components *T, I, F* in general are partially dependent and partially independent.

### 19. Intuitionistic Fuzzy Operators ignore the Indeterminacy, while Neutrosophic Operators give Indeterminacy the same weight as to Truth-Membership and Falsehood-Nonmembership

Indeterminacy in intuitioniostic fuzzy set is ignored by the intuitionistic fuzzy aggregation operators, while the neutrosophic aggregation operators treats the indeterminacy at the same weight as the other two neutrosophic components (truth-membership and falsehood-membership).

Thus, even if we have two single-valued triplets, with the sum of each three components equal to 1 { therefore triplets that may be treated both as *intuitionistic fuzzy triplet*, and *neutrosophic triplet* in the same time (since in neutrosophic environment the sum of the neutrosophic components can be any number between 0 and 3, whence in particular we may take the sum 1) }, after applying the intuitionistic fuzzy aggregation operators we get a different result from that obtained after applying the neutrosophic aggregation operators.

### 20. Intuitionistic Fuzzy Operators and Neutrosophic Operators

Let the intuitioniostic fuzzy operators be denoted as: negation ($\neg_{IF}$), intersection ($\wedge_{IF}$), union ($\vee_{IF}$), and implication ($\rightarrow_{IF}$), and the neutrosophic operators [complement, intersection, union, and implication respectively] be denoted as: negation ($\neg_N$), intersection ($\wedge_N$), union ($\vee_N$), and implication ($\rightarrow_N$).

Let $A_1 = (a_1, b_1, c_1)$ and $A_2 = (a_2, b_2, c_2)$ be two triplets such that $a_1, b_1, c_1, a_2, b_2, c_2 \in [0, 1]$ and $a_1 + b_1 + c_1 = a_2 + b_2 + c_2 = 1$.

The intuitionistic fuzzy operators and neutrosophic operators are based on *fuzzy t-norm* ($\wedge_F$) and *fuzzy t-conorm* ($\vee_F$). We'll take for this article the simplest ones:

$$a_1 \wedge_F a_2 = \min\{a_1, a_2\} \text{ and } a_1 \vee_F a_2 = \max\{a_1, a_2\},$$

where $\wedge_F$ is the fuzzy intersection (t-norm) and $\vee_F$ is the fuzzy union (t-conorm).

For the intuitionistic fuzzy implication and neutrosophic implication, we extend the classical implication:

$$A_1 \rightarrow A_2 \text{ that is classically equivalent to } \neg A_1 \vee A_2,$$

where $\rightarrow$ is the classical implication, $\neg$ the classical negation (complement),

and $\vee$ the classical union,

to the intuitionistic fuzzy environment and respectively to the neutrosophic environment.

But taking other fuzzy t-norm and fuzzy t-conorm, the conclusion will be the same, i.e. the results of intuitionistic fuzzy aggregation operators are different from the results of neutrosophic aggregation operators applied on the same triplets.

*Intuitionistic Fuzzy Aggregation Operators* { the simplest used intuitionistic fuzzy operations }:

Intuitionistic Fuzzy Negation:

$\neg_{IF} (a_1, b_1, c_1) = (c_1, b_1, a_1)$

Intuitionistic Fuzzy Intersection:

$(a_1,b_1,c_1) \wedge_{IF} (a_2,b_2,c_2) = (\min\{a_1,a_2\}, 1-\min\{a_1,a_2\}-\max\{c_1,c_2\}, \max\{c_1,c_2\})$

Intuitionistic Fuzzy Union:

$(a_1,b_1,c_1) \vee_{IF} (a_2,b_2,c_2) = (\max\{a_1,a_2\}, 1-\max\{a_1,a_2\}-\min\{c_1,c_2\}, \min\{c_1,c_2\})$

Intuitionistic Fuzzy Implication:

$(a_1,b_1,c_1) \rightarrow_{IF} (a_2,b_2,c_2)$ is intuitionistically fuzzy equivalent to $\neg_{IF}(a_1,b_1,c_1) \vee_{IF} (a_2,b_2,c_2)$

*Neutrosophic Aggregation Operators* { the simplest used neutrosophic operations }:

Neutrosophic Negation:

$\neg_N (a_1, b_1, c_1) = (c_1, 1-b_1, a_1)$

Neutrosophic Intersection:

$(a_1,b_1,c_1) \wedge_N (a_2,b_2,c_2) = (\min\{a_1,a_2\}, \max\{b_1,b_2\}, \max\{c_1,c_2\})$

Neutrosophic Union:

$(a_1,b_1,c_1) \vee_N (a_2,b_2,c_2) = (\max\{a_1,a_2\}, \min\{b_1,b_2\}, \min\{c_1,c_2\})$

Neutrosophic Implication:

$(a_1,b_1,c_1) \rightarrow_N (a_2,b_2,c_2)$ is neutrosophically equivalent to $\neg_N(a_1,b_1,c_1) \vee_N (a_2,b_2,c_2)$

### 21. Numerical Example of Triplet Components whose Summation is 1

Let $A_1 = (0.3, 0.6, 0.1)$ and $A_2 = (0.4, 0.1, 0.5)$ be two triplets, each having the sum:

$0.3 + 0.6 + 0.1 = 0.4 + 0.1 + 0.5 = 1.$

Therefore, they can both be treated as neutrosophic triplets and as intuitionistic fuzzy triplets simultaneously. We apply both, the intuitionistic fuzzy operators and then the neutrosophic operators and we prove that we get different results, especially with respect with Indeterminacy component that is ignored by the intuitionistic fuzzy operators.

### 21.1. Complement/Negation

Intuitionistic Fuzzy:

$\neg_{IF}(0.3, 0.6, 0.1) = (0.1, 0.6, 0.3)$,

and $\neg_{IF}(0.4, 0.1, 0.5) = (0.5, 0.1, 0.4)$.

Neutrosophic:

$\neg_N(0.3, 0.6, 0.1) = (0.1, 1 - 0.6, 0.3) = (0.1, 0.4, 0.3) \neq (0.1, 0.6, 0.3)$,

and $\neg_N(0.4, 0.1, 0.5) = (0.5, 1 - 0.1, 0.4) = (0.5, 0.9, 0.4) \neq (0.5, 0.1, 0.4)$.

### 21.2. Intersection:

Intuitionistic Fuzzy

$(0.3, 0.6, 0.1) \wedge_{IF} (0.4, 0.1, 0.5) = (\min\{0.3, 0.4\}, 1 - \min\{0.3, 0.4\} - \max\{0.1, 0.5\}, \max\{0.1, 0.5\}) = (0.3, 0.2, 0.5)$

As we see, the indeterminacies 0.6 of $A_1$ and 0.1 of $A_2$ were completely ignored into the above calculations, which is unfair. Herein, the resulting indeterminacy from intersection is just what is left from truth-membership and falsehood-nonmembership { 1 - 0.3 - 0.5 = 0.2 }.

Neutrosophic

$(0.3, 0.6, 0.1) \wedge_N (0.4, 0.1, 0.5) = (\min\{0.3, 0.4\}, \max\{0.6, 0.1\}, \max\{0.1, 0.5\}) = (0.3, 0.6, 0.5) \neq (0.3, 0.2, 0.5)$

In the neutrosophic environment the indeterminacies 0.6 of $A_1$ and 0.1 of $A_2$ are given full consideration in calculating the resulting intersection's indeterminacy: $\max\{0.6, 0.1\} = 0.6$.

### 21.3. Union:

Intuitionistic Fuzzy:

$(0.3, 0.6, 0.1) \vee_{IF} (0.4, 0.1, 0.5) = (\max\{0.3, 0.4\}, 1 - \max\{0.3, 0.4\} - \min\{0.1, 0.5\}, \max\{0.1, 0.5\}) = (0.4, 0.5, 0.1)$

Again, the indeterminacies 0.6 of $A_1$ and 0.1 of $A_2$ were completely ignored into the above calculations, which is not fair. Herein, the resulting indeterminacy from the union is just what is left from truth-membership and falsehood-nonmembership { 1 - 0.4 - 0.1 = 0.5 }.

Neutrosophic:

$(0.3, 0.6, 0.1) \vee_N (0.4, 0.1, 0.5) = (\max\{0.3, 0.4\}, \min\{0.6, 0.1\}, \min\{0.1, 0.5\}) = (0.4, 0.1, 0.1) \neq (0.4, 0.5, 0.1)$

Similarly, in the neutrosophic environment the indeterminacies 0.6 of $A_1$ and 0.1 of $A_2$ are given full consideration in calculating the resulting union's indeterminacy: $\min\{0.6, 0.1\} = 0.1$.

### 21.4. Implication:

Intuitionistic Fuzzy

$(0.3, 0.6, 0.1) \rightarrow_{IF} (0.4, 0.1, 0.5) = \neg_{IF}(0.3, 0.6, 0.1) \vee_{IF} (0.4, 0.1, 0.5) = (0.1, 0.6, 0.3) \vee_{IF} (0.4, 0.1, 0.5) = (0.4, 0.3, 0.3)$

Similarly, indeterminacies of $A_1$ and $A_2$ are completely ignored.

Neutrosophic

$(0.3, 0.6, 0.1) \rightarrow_N (0.4, 0.1, 0.5) = \neg_N(0.3, 0.6, 0.1) \vee_N (0.4, 0.1, 0.5) = (0.1, 0.4, 0.3) \vee_N (0.4, 0.1, 0.5) = (0.4, 0.1, 0.3) \neq (0.4, 0.3, 0.3)$

While in the neutrosophic environment the indeterminacies of $A_1$ and $A_2$ are taken into calculations.

### 21.5. Remark:

We have proven that even when the sum of the triplet components is equal to 1, as demanded by intuitionistic fuzzy environment, the results of the intuitionistic fuzzy operators are different from those of the neutrosophic operators – because the indeterminacy is ignored into the intuitionistic fuzzy operators.

## 22. Simple Counterexample 1, Showing Different Results between Neutrosophic Operators and Intuitionistic Fuzzy Operators Applied on the Same Sets (with component sums > 1 or < 1)

Let the universe of discourse $\mathcal{U} = \{x_1, x_2\}$, and two neutrosophic sets included in $\mathcal{U}$:

$A_N = \{x_1(0.8, 0.3, 0.5), x_2(0.9, 0.2, 0.6)\}$, and

$B_N = \{x_1(0.2, 0.1, 0.3), x_2(0.6, 0.2, 0.1)\}$.

Whence, for $A_N$ one has, after using Atanassov and Vassiliev's transformations (48)' - (51)':

$T_A^{IIFS}(x_1) = \dfrac{0.8}{0.9 + 0.3 + 0.6} = \dfrac{0.8}{1.8} \approx 0.44$;

$I_A^{IIFS}(x_1) = \dfrac{0.3}{1.8} \approx 0.17$;

$F_A^{IIFS}(x_1) = \dfrac{0.5}{1.8} \approx 0.28$.

The refusal degree for $x_1$ with respect to $A_N$ is:

$R_A^{IIFS}(x_1) = 1 - 0.44 - 0.17 - 0.28 = 0.11$.

Then:

$$T_A^{IIFS}(x_2) = \frac{0.9}{1.8} = 0.50;$$

$$I_A^{IIFS}(x_2) = \frac{0.2}{1.8} \approx 0.11;$$

$$F_A^{IIFS}(x_2) = \frac{0.6}{1.8} \approx 0.33.$$

The refusal degree for $x_2$ with respect to $A_N$ is:

$$R_A^{IIFS}(x_2) = 1 - 0.50 - 0.11 - 0.33 = 0.06.$$

Then:

$$A_{IIFS} = \{x_1(0.44, 0.17, 0.28), x_2(0.50, 0.11, 0.33)\}.$$

For $B_N$ one has:

$$T_B^{IIFS}(x_1) = \mu_B^i(x_1) = \frac{0.2}{0.6+0.2+0.3} = \frac{0.2}{1.1} \approx 0.18;$$

$$I_B^{IIFS}(x_1) = v_B^i(x_1) = \frac{0.1}{1.1} \approx 0.09;$$

$$F_B^{IIFS}(x_1) = \pi_B^i(x_1) = \frac{0.3}{1.1} \approx 0.27.$$

The refusal degree for $x_1$ with respect to $B_N$ is:

$$R_B^{IIFS}(x_1) = 1 - 0.18 - 0.09 - 0.27 = 0.46.$$

$$T_B^{IIFS}(x_2) = \frac{0.6}{1.1} \approx 0.55;$$

$$I_B^{IIFS}(x_2) = \frac{0.2}{1.1} \approx 0.18;$$

$$F_B^{IIFS}(x_2) = \frac{0.1}{1.1} \approx 0.09.$$

The refusal degree for $x_2$ with respect to the set $B_N$ is:

$$R_B^{IIFS}(x_2) = 1 - 0.55 - 0.18 - 0.09 = 0.18.$$

Therefore:

$$B_{IIFS} = \{x_1, (0.18, 0.09, 0.27), x_2(0.55, 0.18, 0.09)\}.$$

Therefore, the neutrosophic sets:

$$A_N = \{x_1(0.8, 0.3, 0.5), x_2(0.9, 0.2, 0.6)\} \text{ and}$$

$$B_N = \{x_1(0.2, 0.1, 0.3), x_2(0.6, 0.2, 0.1)\},$$

where transformed (restricted), using Atanassov and Vassiliev's transformations (48)-(51), into inconsistent intuitionistic fuzzy sets respectively as follows:

$$A_{IIFS}^{(t)} = \{x_1(0.44, 0.17, 0.28), x_2(0.50, 0.11, 0.33)\} \text{ and}$$

$$B_{IIFS}^{(t)} = \{x_1(0.18, 0.09, 0.27), x_2(0.55, 0.18, 0.09)\},$$

where the upper script (t) means "after Atanassov and Vassiliev's transformations".

We shall remark that the set $B_N$, as neutrosophic set (where the sum of the components is allowed to also be strictly less than 1 as well), happens to be in the same time an inconsistent intuitionistic fuzzy set, or $B_N \equiv B_{IIFS}$.

Therefore, $B_N$ transformed into $B_{IIFS}^{(t)}$ was a distortion of $B_N$, since we got different IIFS components:

$$x_1^{B_N}(0.2, 0.1, 0.3) \equiv x_1^{B_{IIFS}}(0.2, 0.1, 0.3) \neq x_1^{B_{IIFS}^{(t)}}(0.18, 0.09, 0.27).$$

Similarly:

$$x_2^{B_N}(0.6, 0.2, 0.1) \equiv x_2^{B_{IIFS}}(0.6, 0.2, 0.1) \neq x_2^{B_{IIFS}^{(t)}}(0.55, 0.18, 0.09).$$

Further on, we show that the NS operators and IIFS operators, applied on these sets, give *different results*. For each individual set operation (intersection, union, complement/negation, inclusion/implication, and equality/equivalence) there exist classes of operators, not a single one. We choose the simplest one in each case, which is based on min / max (fuzzy t-norm / fuzzy t-conorm).

### 22.1. Intersection

*Neutrosophic Sets* ( min / max / max )

$$x_1^A \wedge_N x_1^B = (0.8, 0.3, 0.5) \wedge_N (0.2, 0.1, 0.3)$$
$$= (\min\{0.8, 0.2\}, \max\{0.3, 0.1\}, \max\{0.5, 0.3\}) = (0.2, 0.3, 0.5).$$

$$x_2^A \wedge_N x_2^B = (0.9, 0.2, 0.6) \wedge_N (0.6, 0.2, 0.1) = (0.6, 0.2, 0.6).$$

Therefore:

$$A_N \wedge_N B_N = \{x_1(0.2, 0.3, 0.5), x_2(0.6, 0.2, 0.6)\} \stackrel{\text{def}}{=} C_N.$$

*Inconsistent Intuitionistic Fuzzy Set* ( min / max / max )

$$x_1^A \wedge_{IIFS} x_1^B = (0.44, 0.17, 0.28) \wedge_{IIFS} (0.18, 0.09, 0.27) =$$
$$(\min\{0.44, 0.18\}, \max\{0.17, 0.09\}, \max\{0.28, 0.27\}) =$$
$$(0.18, 0.17, 0.28)$$

$$x_2^A \wedge_{IIFS} x_2^B = (0.50, 0.11, 0.33) \wedge_{IIFS} (0.55, 0.18, 0.09)$$
$$= (0.50, 0.18, 0.33).$$

Since in IIFS the sum of components is not allowed to surpass 1, we normalize:

$$\left(\frac{0.50}{1.01}, \frac{0.11}{1.01}, \frac{0.33}{1.01}\right) \approx (0.495, 0.109, 0.326).$$

Therefore:

$$A_{IIFS} \wedge_{IIFS} B_{IIFS} = \{x_1(0.18, 0.17, 0.28), x_2(0.495, 0.109, 0.326)\} \stackrel{\text{def}}{=} C_{IIFS}.$$

Also:

$$T_{A_N \wedge_N B_N}(x_1) = 0.2 < 0.3 = I_{A_N \wedge_N B_N}(x_1),$$

while

$$T_{A_{IIFS} \wedge_{IIFS} B_{IIFS}}(x_1) = 0.18 > 0.17 = I_{A_{IIFS} \wedge_{IIFS} B_{IIFS}}(x_1),$$

and other discrepancies can be seen.

*Inconsistent Intuitionistic Fuzzy Set* ( with min / min / max, as used by Cuong [20] in order to avoid the sum of components surpassing 1; but this is in discrepancy with the IIFS/PFS union that uses max / min / min, not max / max / min ):

$$x_1^A \wedge_{IIFS2} x_1^B = (0.44, 0.17, 0.28) \wedge_{IIFS2} (0.18, 0.09, 0.27) =$$
$$(\min\{0.44, 0.18\}, \min\{0.17, 0.09\}, \max\{0.28, 0.27\}) =$$
$$(0.18, 0.09, 0.28)$$

$$x_2^A \wedge_{IIFS2} x_2^B = (0.50, 0.11, 0.33) \wedge_{IIFS2} (0.55, 0.18, 0.09)$$
$$= (0.50, 0.11, 0.33).$$

Therefore:

$$A_{IIFS} \wedge_{IIFS2} B_{IIFS} = \{x_1(0.18, 0.09, 0.28), x_2(0.50, 0.11, 0.33)\} \stackrel{\text{def}}{=} C_{IIFS2}$$

We see that:

$$A_N \wedge_N B_N \neq A_{IIFS} \wedge_{IIFS} B_{IIFS}, \text{ or } C_N \neq C_{IIFS};$$

and $A_N \wedge_N B_N \neq A_{IIFS} \wedge_{IIFS2} B_{IIFS}$, $C_N \neq C_{IIFS2}$. Also $C_{IIFS} \neq C_{IIFS2}$.

Let's transform the above neutrosophic set $C_N$, resulted from the application of the neutrosophic intersection operator,

$$C_N = \{x_1(0.2, 0.3, 0.5), x_2(0.6, 0.2, 0.6)\},$$

into an inconsistent intuitionistic fuzzy set, employing the same equations (48) – (50) of transformations [denoted by (t)], provided by Atanassov and Vassiliev, which are equivalent {using (T, I, F)-notations} to (48)'-(50)'

$$(t)T_C^{IIFS}(x_1) = \frac{0.2}{0.6+0.3+0.6} = \frac{0.2}{1.5} \simeq 0.13;$$

$$(t)I_C^{IIFS}(x_1) = \frac{0.3}{1.5} = 0.20;$$

$$(t)F_C^{IIFS}(x_1) = \frac{0.5}{1.5} \simeq 0.33.$$

$$(t)T_C^{IIFS}(x_2) = \frac{0.6}{1.5} \simeq 0.40;$$

$$(t)I_C^{IIFS}(x_2) = \frac{0.2}{1.5} \simeq 0.13;$$

$$(t)F_C^{IIFS}(x_2) = \frac{0.6}{1.5} \simeq 0.40.$$

Whence the results of neutrosophic and IIFS/PFS are totally different:

$$C_{IIFS}^{(t)} = \{x_1(0.13, 0.20, 0.33), x_2(0.40, 0.13, 0.40)\} \neq$$

$$\{x_1(0.18, 0.17, 0.28), x_2(0.495, 0.109, 0.326)\} \equiv C_{IIFS},$$

and

$$C_{IIFS}^{(t)} \neq \{x_1(0.18, 0.09, 0.28), x_2(0.50, 0.11, 0.33)\} = C_{IIFS2}.$$

22.2. Union

*Neutrosophic Sets ( max / min / min )*

$$x_1^A \vee_N x_1^B = (0.8, 0.3, 0.5) \vee_N (0.2, 0.1, 0.3)$$
$$= (\max\{0.8, 0.2\}, \min\{0.3, 0.1\}, \min\{0.5, 0.3\}) = (0.8, 0.1, 0.3).$$

$$x_2^A \vee_N x_2^B = (0.9, 0.2, 0.6) \vee_N (0.6, 0.2, 0.1) = (0.9, 0.2, 0.1).$$

Therefore:

$$A_N \vee_N B_N = \{x_1(0.8, 0.1, 0.3), x_2(0.9, 0.2, 0.1)\} \stackrel{\text{def}}{=} D_N.$$

*Inconsistent Intuitionistic Fuzzy Sets ( max / min / min [3] )*

$$x_1^A \vee_{IIFS} x_1^B = (0.44, 0.17, 0.28) \vee_{IIFS} (0.18, 0.09, 0.27)$$
$$= (\max\{0.44, 0.18\}, \min\{0.17, 0.09\}, \min\{0.28, 0.27\})$$
$$= (0.44, 0.09, 0.27).$$

$$x_2^A \vee_{IIFS} x_2^B = (0.50, 0.11, 0.33) \vee_{IIFS} (0.55, 0.18, 0.09)$$
$$= (0.55, 0.11, 0.09).$$

Therefore:

$$A_{IIFS} \vee_{IIFS} B_{IIFS} = \{x_1(0.44, 0.09, 0.27), x_2(0.55, 0.11, 0.09)\} \stackrel{\text{def}}{=} D_{IIFS}$$

a) We see that the results are totally different:

$$A_N \vee_N B_N \neq A_{IIFS} \vee_{IIFS} B_{IIFS}, \text{ or } D_N \neq D_{IIFS}.$$

b) Let's transform the above neutrosophic set, $D_N$, resulted from the application of neutrosophic union operator,

$$D_N = \{x_1(0.8, 0.1, 0.3), x_2(0.9, 0.2, 0.1)\},$$

into an inconsistent intuitionistic fuzzy set, employing the same equations (48) -(50) of transformation [ denoted by (*t*) ], provided by Atanassov and Vassiliev, which are equivalent [using (T, I, F) notations] to (48)'-(50)':

$$(t)T_D^{IIFS}(x_1) = \frac{0.8}{0.9+0.2+0.3} = \frac{0.8}{1.4} \simeq 0.57;$$

$$(t)I_D^{IIFS}(x_1) = \frac{0.1}{1.4} \simeq 0.07;$$

$$(t)I_D^{IIFS}(x_1) = \frac{0.3}{1.4} \simeq 0.21.$$

$$(t)T_D^{IIFS}(x_2) = \frac{0.9}{1.4} \simeq 0.64;$$

$$(t)I_D^{IIFS}(x_2) = \frac{0.2}{1.4} \simeq 0.14;$$

$$(t)F_D^{IIFS}(x_2) = \frac{0.1}{1.4} \simeq 0.07.$$

Whence:

$$D_{IIFS}^{(t)} = \{x_1(0.57, 0.07, 0.21), x_2(0.64, 0.14, 0.07)\}$$
$$\neq \{x_1(0.44, 0.09, 0.27), x_2(0.55, 0.11, 0.09)\} \equiv D_{IIFS}.$$

The results again are totally different.

### 22.3. Corollary

Therefore, no matter if we first transform the neutrosophic components into inconsistent intuitionistic fuzzy components (as suggested by Atanassov and Vassiliev) and then apply the IIFS operators, or we first apply the neutrosophic operators on neutrosophic components, and then later transform the result into IIFS components, in both ways the obtained results in the neutrosophic environment are totally different from the results obtained in the IIFS environment.

### 24. Normalization

Further on, the authors propose the normalization of the neutrosophic components, where Atanassov and Vassiliev's [6] equations (57) – (59) are equivalent, using neutrosophic notations, to the following.

Let $\mathcal{U}$ be a universe of discourse, a set $A \subseteq \mathcal{U}$, and a generic element $x \in \mathcal{U}$, with the neutrosophic components:

$x(T_A(x), I_A(x), F_A(x))$, where

$T_A(x), I_A(x), F_A(x) \in [0, 1]$, and

$T_A(x) + I_A(x) + F_A(x) \leq 3$, for all $x \in U$.

Suppose $T_A(x) + I_A(x) + F_A(x) \neq 0$, for all $x \in U$.

Then, by the below normalization of neutrosophic components, Atanassov and Vassiliev obtain the following intuitionistic fuzzy components $(T_A^{IFS}, I_A^{IFS}, F_A^{IFS})$:

$$T_A^{IFS}(x) = \frac{T_A(x)}{T_A(x) + I_A(x) + F_A(x)} \in [0, 1] \tag{57'}$$

$$I_A^{IFS}(x) = \frac{I_A(x)}{T_A(x) + I_A(x) + F_A(x)} \in [0, 1] \tag{58'}$$

$$F_A^{IFS}(x) = \frac{F_A(x)}{T_A(x) + I_A(x) + F_A(x)} \in [0, 1] \tag{59'}$$

and

$T_A^{IFS}(x) + I_A^{IFS}(x) + F_A^{IFS}(x) = 1$, for all $x \in U$.

### 16.1. Counterexample 2

Let's come back to the previous *Counterexample 1*.

$U = \{x_1, x_2\}$ be a universe of discourse, and let two neutrosophic sets included in $U$:

$A_N = \{x_1(0.8, 0.3, 0.5), x_2(0.9, 0.2, 0.6)\}$, and

$B_N = \{x_1(0.2, 0.1, 0.3), x_2(0.6, 0.2, 0.1)\}$.

Let's normalize their neutrosophic components, as proposed by Atanassov and Vassiliev, in order to restrain them to intuitionistic fuzzy components:

$$A_{IFS} = \left\{ x_1 \left( \frac{0.8}{0.8 + 0.3 + 0.5}, \frac{0.3}{1.6}, \frac{0.5}{1.6} \right), x_2 \left( \frac{0.9}{1.7}, \frac{0.2}{1.7}, \frac{0.6}{1.7} \right) \right\}$$

$$\approx \{x_1(0.50, 0.19, 0.31), x_2(0.53, 0.12, 0.35)\}$$

$$\equiv \{x_1(0.50, 0.31), x_2(0.53, 0.35)\},$$

since the indeterminacy (called *hesitant degree* in IFS) is neglected.

$$B_{IFS} = \left\{ x_1 \left( \frac{0.2}{0.6}, \frac{0.1}{0.6}, \frac{0.3}{0.6} \right), x_2 \left( \frac{0.6}{0.9}, \frac{0.2}{0.9}, \frac{0.1}{0.9} \right) \right\}$$

$$\approx \{x_1(0.33, 0.17, 0.50), x_2(0.67, 0.22, 0.11)\}$$

$$\equiv \{x_1(0.33, 0.50), x_2(0.67, 0.11)\},$$

since the indeterminacy (hesitance degree) is again neglected.

The intuitionistic fuzzy operators are applied only on truth-membership and false-nonmembership (but not on indeterminacy).

### 24.1.1. Intersection

*Intuitionistic Fuzzy Intersection ( min / max )*

$$x_1^A \wedge_{IFS} x_1^B = (0.50, 0.31) \wedge_{IFS} (0.33, 0.50) =$$

$$(min\{0.50, 0.33\}, max\{0.31, 0.50\}) = (0.33, 0.50) \equiv (0.33, 0.17, 0.50)$$

,

after adding the indeterminacy which is what's left up to 1, i.e. $1 - 0.33 - 0.50 = 0.17$.

$$x_2^A \wedge_{IFS} x_2^B = (0.53, 0.35) \wedge_{IFS} (0.67, 0.11)$$

$$= (min\{0.53, 0.63\}, max\{0.35, 0.11\}) = (0.53, 0.35)$$

$$\equiv (0.53, 0.12, 0.35),$$

after adding the indeterminacy.

The results of NS and IFS intersections are clearly very different:

$$A_N \wedge_N B_N = \{x_1(0.2, 0.3, 0.5), x_2(0.6, 0.2, 0.6)\}$$

$$\neq \{x_1(0.33, 0.17, 0.50), x_2(0.53, 0.12, 0.35)\} = A_{IFS} \wedge_{IFS} B_{IFS}$$

Even more distinction, between the NS intersection and IFS intersection of the same elements (whose sums of components equal 1) $x_1^A = (0.50, 0.19, 0.31)$ and $x_1^B = (0.33, 0.17, 0.50)$ one obtains unequal results, using the (min / max / max) operator:

$$x_1^A \wedge_N x_1^B = (0.50, 0.19, 0.31) \wedge_N (0.33, 0.17, 0.50) = (0.33, 0.19, 0.50),$$

while

$$x_1^A \wedge_{IFS} x_1^B = (0.50, 0.19, 0.31) \wedge_{IFS} (0.33, 0.17, 0.50)$$

$$\equiv (0.50, 0.31) \wedge_{IFS} (0.33, 0.50) \text{ \{after ignoring the indeterminacy in IFS\}}$$

$$= (0.33, 0.50) \equiv (0.33, 0.17, 0.50) \neq (0.33, 0.19, 0.50).$$

*24.1.2. Union*

*Intuitionistic Fuzzy Union ( max / min / min )*

$$x_1^A \vee_{IFS} x_1^B = (0.50, 0.31) \vee_{IFS} (0.33, 0.50)$$
$$= (\max\{0.50, 0.33\}, \min\{0.31, 0.50\}) = (0.50, 0.31)$$
$$\equiv (0.50, 0.19, 0.31),$$

after adding the indeterminacy.

$$x_2^A \vee_{IFS} x_2^B = (0.53, 0.35) \vee_{IF} (0.67, 0.11)$$
$$= (\max\{0.53, 0.67\}, \min\{0.35, 0.11\}) = (0.67, 0.11)$$
$$\equiv (0.67, 0.22, 0.11),$$

after adding the indeterminacy.

The results of NS and IFS unions are clearly very different:

$$A_N \vee_N B_N = \{x_1(0.8, 0.1, 0.3), x_2(0.9, 0.2, 0.1)\}$$
$$\neq \{x_1(0.50, 0.19, 0.31), x_2(0.67, 0.22, 0.11)\} = A_{IFS} \vee_{IF} B_{IFS}.$$

Even more distinction, for the NS and IFS union of the same elements:

$$x_1^A \vee_N x_1^B = (0.50, 0.19, 0.31) \vee_N (0.33, 0.17, 0.50) = (0.50, 0.17, 0.31),$$

while

$$x_1^A \vee_{IFS} x_1^B = (0.50, 0.19, 0.31) \vee_{IFS} (0.33, 0.17, 0.50) \equiv$$
$$(0.50, 0.31) \vee_{IFS} (0.33, 0.50)$$
$$=$$
$$(0.50, 0.31) \equiv (0.50, 0.19, 0.31) \{\text{after adding indeterminacy}\} \neq (0.50, 0.17, 0.31).$$

**25. Indeterminacy Makes a Big Difference between NS and IFS**

The authors [6] assert that,

"Therefore, the NFS can be *represented* by an IFS" (page 5),

but this is not correct, since it should be:

The NFS (neutrosophic fuzzy set ≡ single-valued neutrosophic set) can be *restrained* (degraded) to an IFS (intuitionistic fuzzy set), yet the <u>independence of components</u> is lost and the results of the aggregation operators are totally different between the neutrosophic environment and intuitionistic fuzzy environment, since Indeterminacy is ignored by IFS operators.

Since in single-valued neutrosophic set the neutrosophic components are *independent* (their sum can be up to 3, and if a component increases or decreases, it does not change the others), while in intuitionistic fuzzy set the components are *dependent* (in general if one changes, one or both the other components change in order to keep their sum equal to 1). Also, applying the neutrosophic operators is a better aggregation since the indeterminacy (*I*) is involved into all neutrosophic (complement/negation, intersection, union, inclusion/inequality/implication, equality/equivalence) operators while all intuitionistic fuzzy operators *ignore* (do not take into calculation) the indeterminacy.

That is why the results after applying the neutrosophic operators and intuitionistic fuzzy operators on the same sets are *different* as proven above.

### 26. Paradoxes cannot be Represented by the Intuitionistic Fuzzy Logic

No previous set/logic theories, including IFS or Intuitionistic Fuzzy Logic (IFL), since the sum of components was not allowed above 1, could characterize a paradox, which is a proposition that is true (T = 1) and false (F = 1) simultaneously, therefore the paradox is 100% indeterminate (I = 1). In Neutrosophic Logic (NL) a paradoxical proposition $P_{NL}$ is represented as: $P_{NL}(1, 1, 1)$.

If one uses Atanassov and Vassiliev's transformations (for example the normalization) [6], we get $P_{IFL}(1/3, 1/3, 1/3)$, but this one cannot represent a paradox, since a paradox is 100% true and 100% false, not 33% true and 33% false.

### 27. Single-Valued Atanassov's Intuitionistic Fuzzy Set of second type, also called Single-Valued Pythagorean Fuzzy Set

Single-Valued Atanassov's Intuitionistic Fuzzy Sets of second type (AIFS2) [23], also called Single-Valued Pythagorean Fuzzy Set (PyFS) [24], is defined as follows (using T, I, F notations for the components):

Definition of IFS2 (PyFS)

It is a set $A_{AIFS2} \equiv A_{PyFS}$ from the universe of discourse U such that:

$A_{AIFS2} \equiv A_{PyFS} = \{<x, T_A(x), F_A(x)> | x \in U\}$,

where, for all $x \in U$, the functions $T_A(x), F_A(x) : U \to [0, 1]$, represent the degree of membership (truth) and degree on nonmembership (falsity) respectively, that satisfy the conditions:

$$0 \leq T_A^2(x) + F_A^2(x) \leq 1,$$

whence the hesitancy degree is:

$$I_A(x) = \sqrt{1 - T_A^2(x) - F_A^2(x)} \in [0,1].$$

### 28. Single-Valued Refined Pythagorean Fuzzy Set (RPyFS)

We propose now for the first time the Single-Valued Refined Pythagorean Fuzzy Set (RPyFS):

$$A_{RAIFS2} = A_{RPyFS} = \{x(T_A^1(x), T_A^2(x), ..., T_A^p(x); F_A^1(x), F_A^2(x), ..., F_A^s(x)), p+s \geq 3, x \in U\}$$

where *p* and *s* are positive nonzero integers, and for all $x \in U$, the functions $T_A^1(x), T_A^2(x), ..., T_A^p(x), F_A^1(x), F_A^2(x), ..., F_A^s(x) : U \to [0, 1]$, represent the degrees of sub-membership (sub-truth) of types *1, 2, ..., p,* and degrees on sub-nonmembership (sub-falsity) of types *1, 2, ..., s* respectively, that satisfy the condition:

$$0 \leq \sum_{1}^{p}(T_A^j)^2 + \sum_{1}^{s}(F_A^l)^2 \leq 1,$$

whence the refined hesitancy degree is:

$$I_A(x) = \sqrt{1 - \sum_{1}^{p}(T_A^j)^2 - \sum_{1}^{s}(F_A^l)^2} \in [0,1].$$

The Single-Valued Refined Pythagorean Fuzzy Set is a particular case of the Single-Valued Refined Neutrosophic Set.

### 29. The components of Atanassov's Intuitionistic Fuzzy Set of second type (Pythagorean Fuzzy Set) are not Independent

Princy R and Mohana K assert in [23] that:

> "the truth and falsity values and hesitancy value can be independently considered as membership and non-membership and hesitancy degrees respectively".

But this is untrue, since in IFS2 (PyFS) the components are not independent, because they are connected (dependent on each other) through this inequality:

$$T_A^2(x) + F_A^2(x) \leq 1.$$

### 30. Let's see a **Counterexample 3**:

If T = 0.9, then $T^2 = 0.9^2 = 0.81$, whence $F^2 \leq 1 - T^2 = 1 - 0.81 = 0.19$,
or $F \leq \sqrt{0.19} \approx 0.44$.

Therefore, if T = 0.9, then F is restricted to be less than equal to $\sqrt{0.19}$.

While in NS if T = 0.9, F can be equal to any number in [0, 1], F can be even equal to 1.

Also, hesitancy degree clearly depends on T and F, because the formula of hesitancy degree is an equation depending on T and F, as below:

$$I_A(x) = \sqrt{1 - T_A^2(x) - F_A^2(x)} \in [0,1].$$

If T = 0.9 and F = 0.2, then hesitancy

$$I = \sqrt{1 - 0.9^2 - 0.2^2} = \sqrt{0.15} \approx 0.39.$$

Again, in NS if T = 0.9 and F = 0.2, I can be equal to any number in [0, 1], not only to $\sqrt{0.15}$.

### 31. Neutrosophic Set is a Generalization of Pythagorean Fuzzy Set

In the definition of PyFS, one has $T_A(x), F_A(x) \in [0, 1]$, which involves that $T_A(x)^2, F_A(x)^2 \in [0, 1]$ too;

we denote $T_A^{NS}(x) = T_A(x)^2, F_A^{NS}(x) = F_A(x)^2$, and

$I_A^{NS}(x) = I_A(x)^2 = 1 - T_A(x)^2 - F_A(x)^2 \in [0,1]$, where "NS" stands for Neutrosophic Set.

Therefore, one gets: $T_A^{NS}(x) + I_A^{NS}(x) + F_A^{NS}(x) = 1$,

which is a particular case of the neutrosophic set, since in NS the sum of the components can be any number between 0 and 3, hence into PyFS has been chosen the sum of the components be equal to 1.

### 32. Spherical Fuzzy Set (SFS)
### Definition of Spherical Fuzzy Set

A Single-Valued Spherical Fuzzy Set (SFS) [25, 26], of the universe of discourse U, is defined as follows:

$A_{SFS} = \{<x, T_A(x), I_A(x), F_A(x)> | x \in U\}$,

where, for all x ∈ U, the functions $T_A(x)$, $I_A(x)$, $F_A(x) : U \to [0, 1]$, represent the degree of membership (truth), the degree of hesitancy, and degree on nonmembership (falsity) respectively, that satisfy the conditions:

$$0 \leq T_A^2(x) + I_A^2(x) + F_A^2(x) \leq 1,$$

whence the refusal degree is:

$$R_A(x) = \sqrt{1 - T_A^2(x) - I_A^2(x) - F_A^2(x)} \in [0,1].$$

### 33. Single-Valued n-HyperSpherical Fuzzy Set (n-HSFS)

Smarandache (2019) generalized for the first time the spherical fuzzy set to n-hyperspherical fuzzy set.

*Definition of n-HyperSpherical Fuzzy Set.*

A Single-Valued n-HyperSpherical Fuzzy Set (n-HSFS), of the universe of discourse U, is defined as follows:

$$A_{n\text{-HSFS}} = \{<x, T_A(x), I_A(x), F_A(x)> \mid x \in U\},$$

where, for all x ∈ U, the functions $T_A(x)$, $I_A(x)$, $F_A(x) : U \to [0, 1]$, represent the degree of membership (truth), the degree of hesitancy, and degree on nonmembership (falsity) respectively, that satisfy the conditions:

$$0 \leq T_A^n(x) + I_A^n(x) + F_A^n(x) \leq 1, \text{ for } n \geq 1,$$

whence the refusal degree is:

$$R_A(x) = \sqrt{1 - T_A^n(x) - I_A^n(x) - F_A^n(x)} \in [0,1].$$

It is clear that 2-HyperSpherical Fuzzy Set (i.e. when *n = 2*) is a spherical fuzzy set.

### 34. The n-HyperSpherical Fuzzy Set is a particular case of the Neutrosophic Set.

Because, $T_A(x), I_A(x), F_A(x) \in [0, 1]$ implies that, for *n ≥ 1* one has $T_A^n(x), I_A^n(x), F_A^n(x) \in [0,1]$ too, so they are neutrosophic components as well; therefore each n-HSFS is a NS.

But the reciprocal is not true, since if at least one component is 1 and from the other two components at least one is > 0, for example $T_A(x) = 1$, and $I_A(x) > 0$, $F_A(x) \in [0, 1]$, then $T_A^n(x) + I_A^n(x) + F_A^n(x) > 1$ for *n ≥ 1*. Therefore, there are infinitely many triplets T, I, F that are NS components, but they are not n-HSFS components.

### 35. The components of the Spherical Fuzzy Set are not Independent

Princy R and Mohana K assert in [23] that:

> "In spherical fuzzy sets, while the squared sum of membership, non-membership and hesitancy parameters can be between 0 and 1, each of them can be defined between 0 and 1 independently."

But this is again untrue, the above parameters cannot be defined independently.

### 36. Counterexample 4

If *T = 0.9* then *F* cannot be for example equal to *0.8,*
since $0.9^2 + 0.8^2 = 1.45 > 1$,
but the sum of the squares of components is not allowed to be greater than 1.
So *F* depends on *T* in this example.

Two components are independent if no matter what value gets one component will not affect the other component's value.

### 37. Neutrosophic Set is a generalization of the Spherical Fuzzy Set

In [25] Gündoğlu and Kahraman assert about:

> "superiority of SFS [i.e. Spherical Fuzzy Set] with respect to Pythagorean, intuitionistic fuzzy and neutrosophic sets";
> also:
> "SFSs are a generalization of Pythagorean Fuzzy Sets (PFS) and neutrosophic sets".

While it is *true* that the spherical fuzzy set is a generalizations of Pythagorean fuzzy set and of intuitionistic fuzzy set, it is *false* that spherical fuzzy set is a generalization of neutrosophic set.
Actually it's the opposite: neutrosophic set is a generalization of spherical fuzzy set. We prove it bellow.

*Proof*
In the definition of the spherical fuzzy set one has:
$T_A(x), I_A(x), F_A(x) \in [0, 1]$, which involves that $T_A(x)^2, I_A(x)^2, F_A(x)^2 \in [0, 1]$ too.
Let's denote: $T_A^{NS}(x) = T_A(x)^2, I_A^{NS}(x) = I_A(x)^2, F_A^{NS}(x) = F_A(x)^2$, where "NS" stands for neutrosophic set, whence we obtain, using SFS definition:
$0 \leq T_A^{NS}(x) + I_A^{NS}(x) + F_A^{NS}(x) \leq 1$,
which is a particular case of the single-valued neutrosophic set, where the sum of the components T, I, F can be any number between 0 and 3. So now we can choose the sum up to 1.

### 38. Counterexample 5

If we take $T_A(x) = 0.9$, $I_A(x) = 0.4$, $F_A(x) = 0.5$, for some given element *x*, which are neutrosophic components, they are not spherical fuzzy set components because *$0.9^2$ + $0.4^2$ + $0.5^2$ = 1.22 > 1*.
There are infinitely many values for *$T_A(x)$, $I_A(x)$, $F_A(x)$ in [0, 1]* whose sum of squares is strictly greater than 1, therefore they are not spherical fuzzy set components, but they are neutrosophic components.

The elements of a spherical fuzzy set form a 1/8 of a sphere of radius 1, centred into the origin O(0,0,0) of the Cartesian system of coordinates, on the positive *Ox (T), Oy (I), Oz (F)* axes.

While the standard neutrosophic set is a cube of side 1, that has the vertexes: *(0,0,0), (1,0,0), (0,1,0), (0,0,1), (1,1,0), (1,0,1), (0,1,1), (1,1,1)*.

The neutrosophic cube strictly includes the 1/8 fuzzy sphere.

### 39. Single-Valued Refined Spherical Fuzzy Set (RSFS)

We introduce now for the first time the Single-Valued Refined Spherical Fuzzy Set.

$$A_{RSFS} = \{x(T_A^1(x), T_A^2(x),...,T_A^p(x); I_A^1(x), I_A^2(x),...,I_A^r(x);$$
$$F_A^1(x), F_A^2(x),...,F_A^s(x)), p+r+s \geq 4, x \in U\},$$

where *p, r, s* are nonzero positive integers, and for all $x \in U$, the functions
$T_A^1(x), T_A^2(x),...,T_A^p(x), I_A^1(x), I_A^2(x),...,I_A^r(x), F_A^1(x), F_A^2(x),...,F_A^s(x): U \to [0, 1]$,
represent the degrees of sub-membership (sub-truth) of types *1, 2, ..., p*, the degrees of sub-hesitancy of types *1, 2, ..., r*, and degrees on sub-nonmembership (sub-falsity) of types *1, 2, ..., s* respectively, that satisfy the condition:

$$0 \leq \sum_{1}^{p}(T_A^j)^2 + \sum_{1}^{s}(I_A^k)^2 + \sum_{1}^{s}(F_A^l)^2 \leq 1,$$

whence the refined refusal degree is:

$$R_A(x) = \sqrt{1 - \sum_{1}^{p}(T_A^j)^2 - \sum_{1}^{s}(I_A^k)^2 - \sum_{1}^{s}(F_A^l)^2} \in [0,1].$$

The Single-Valued Refined Spherical Fuzzy Set is a particular case of the Single-Valued Refined Neutrosophic Set.

### 40. Single-Valued Spherical Neutrosophic Set

Spherical Neutrosophic Set (SNS) was introduced by Smarandache [27] in 2017.

A Single-Valued Spherical Neutrosophic Set (SNS), of the universe of discourse U, is defined as follows:

$A_{SNS} = \{<x, T_A(x), I_A(x), F_A(x)> | x \in U\}$,

where, for all $x \in U$, the functions $T_A(x), I_A(x), F_A(x) : U \to [0, \sqrt{3}]$, represent the degree of membership (truth), the degree of indeterminacy, and degree on nonmembership (falsity) respectively, that satisfy the conditions:

$0 \leq T_A^2(x) + I_A^2(x) + F_A^2(x) \leq 3$.

The Spherical Neutrosophic Set is a generalization of Spherical Fuzzy Set, because we may restrain the SNS's components to the unit interval $T_A(x), I_A(x), F_A(x) \in [0, 1]$, and the sum of the squared components to 1, i.e. $0 \leq T_A^2(x) + I_A^2(x) + F_A^2(x) \leq 1$.

Further on, if replacing $I_A(x) = 0$ into the Spherical Fuzzy Set, we obtain as particular case the Pythagorean Fuzzy Set.

### 41. Single-Valued n-HyperSpherical Neutrosophic Set (n-HSNS)

Definition of n-HyperSpherical Neutrosophic Set (Smarandache, 2019)

We introduce now for the first time the Single-Valued n-HyperSpherical Neutrosophic Set (n-HSNS), which is a generalization of the Spherical Neutrosophic Set and of n-HyperSpherical Fuzzy Set, of the universe of discourse U, for $n \geq 1$, is defined as follows:

$A_{n-HNS} = \{<x, T_A(x), I_A(x), F_A(x)> | x \in U\}$,

where, for all $x \in U$, the functions $T_A(x), I_A(x), F_A(x) : U \to [0, \sqrt[n]{3}]$, represent the degree of membership (truth), the degree of indeterminacy, and degree on nonmembership (falsity) respectively, that satisfy the conditions:

$0 \leq T_A^n(x) + I_A^n(x) + F_A^n(x) \leq 3$.

## 42. Single-Valued Refined Refined n-HyperSpherical Neutrosophic Set (R-n-HSNS)

We introduce now for the first time the Single-Valued Refined n-HyperSpherical Neutrosophic Set (R-n-HSNS), which is a generalization of the n-HyperSpherical Neutrosophic Set and of Refined n-HyperSpherical Fuzzy Set.

On the universe of discourse U, for $n \geq 1$, we define it as:

$$A_{R-n-HSNS} = \{x(T_A^1(x), T_A^2(x), ..., T_A^p(x); I_A^1(x), I_A^2(x), ..., I_A^r(x);$$
$$F_A^1(x), F_A^2(x), ..., F_A^s(x)), p + r + s \geq 4, x \in U\},$$

where *p, r, s* are nonzero positive integers, and for all $x \in U$, the functions $T_A^1(x), T_A^2(x), ..., T_A^p(x), I_A^1(x), I_A^2(x), ..., I_A^r(x), F_A^1(x), F_A^2(x), ..., F_A^s(x): U \to [0, m^{1/n}]$, represent the degrees of sub-membership (sub-truth) of types *1, 2, ..., p,* the degrees of sub-indeterminacy of types *1, 2, ..., r,* and degrees on sub-nonmembership (sub-falsity) of types *1, 2, ..., s* respectively, that satisfy the condition:

$$0 \leq \sum_{1}^{p}(T_A^j)^n + \sum_{1}^{r}(I_A^k)^n + \sum_{1}^{s}(F_A^l)^n \leq m, \text{ where } p + r + s = m.$$

## 43. Neutrosophic Set is a Generalization of q-Rung Orthopair Fuzzy Set (q-ROFS).

*Definition of q-Rung Orthopair Fuzzy Set.*

Using the same *T, I, F* notations one has as follows.

A *Single-Valued q-Rung Orthopair Fuzzy Set (q-ROFS)* [28], of the universe of discourse U, for a given real number $q \geq 1$, is defined as follows:

$A_{q\text{-ROFS}} = \{<x, T_A(x), F_A(x)> \mid x \in U\}$,

where, for all $x \in U$, the functions $T_A(x), F_A(x) : U \to [0, 1]$, represent the degree of membership (truth), and degree on nonmembership (falsity) respectively, that satisfy the conditions:

$$0 \leq T_A(x)^q + F_A(x)^q \leq 1.$$

Since $T_A(x), F_A(x) \in [0, 1]$, then for any real number $q \geq 1$ one has $T_A(x)^q, F_A(x)^q \in [0,1]$ too.

Let's denote: $T_A^{NS}(x) = T_A(x)^q, F_A^{NS}(x) = F_A(x)^q$, whence it results that:

$0 \leq T_A^{NS}(x) + F_A^{NS}(x) \leq 1$, where what's left may be Indeterminacy.

But this is a particular case of the neutrosophic set, where the sum of components T, I, F can be any number between 0 and 3, and for q-ROFS is it taken to be up to 1. Therefore, any Single-Valued q-Rung Orthopair Fuzzy Set is also a Neutrosophic Set, but the reciprocal is not true. See the next counterexample.

## 44. Counterexample 6.

Let's consider a real number $1 \leq q < \infty$, and a set of single-valued triplets of the form

*(T, I, F)*, with *T, I, F ∈ [0, 1]* that represent the components of the elements of a given set.

The components of the form *(1, F)*, with *F > 0*, and of the form *(T, 1)*, with *T > 0*, constitute NS components as follows: *(1, I, F)*, with *F > 0* and any *I ∈ [0, 1]*, and respectively

*(T, I, 1)*, with *T > 0* and any *I ∈ [0, 1]*, since the sum of the components is allowed to be greater than *1*, i.e. *1 + I + F > 1* and respectively *T + I + 1 > 1*.

But they cannot be components of the elements of a q-ROFS set, since:

$1^q + F^q = 1 + F^q > 1$, because $F > 0$ and $1 \leq q < \infty$; but in q-ROFS the sum has to be ≤ *1*.

Similarly, $T^q + 1^q = T^q + 1 > 1$, because $T > 0$ and $1 \leq q < \infty$; but in q-ROFS the sum has to be ≤ *1*.

### 45. Refined q-Rung Orthopair Fuzzy Set (R-q-ROFS)

We propose now for the first time the Single-Valued Refined q-Rung Orthopair Fuzzy Set (R-q-ROFS):
$$A_{R-q-ROFS} = \{x(T_A^1(x), T_A^2(x), ..., T_A^p(x), F_A^1(x), F_A^2(x), ..., F_A^s(x)), p+s \geq 3, x \in U\},$$
where *p* and *s* are positive nonzero integers, and for all *x ∈ U*, the functions
$T_A^1(x), T_A^2(x), ..., T_A^p(x), F_A^1(x), F_A^2(x), ..., F_A^s(x): U \rightarrow [0, 1]$, represent the degrees of sub-membership (sub-truth) of types *1, 2, ..., p,* and degrees on sub-nonmembership (sub-falsity) of types *1, 2, ..., s* respectively, that satisfy the condition:
$$0 \leq \sum_{1}^{p}(T_A^j)^q + \sum_{1}^{s}(F_A^l)^q \leq 1, \text{ for } q \geq 1,$$
whence the refined hesitancy degree is:
$$I_A(x) = [1 - \sum_{1}^{p}(T_A^j)^q - \sum_{1}^{s}(F_A^l)^q]^{1/q} \in [0,1].$$

The Single-Valued Refined q-Rung Fuzzy Set is a particular case of the Single-Valued Refined Neutrosophic Set.

### 46. Regret Theory is a Neutrosophication Model

Regret Theory (2010) [29] is actually a Neutrosophication (1998) Model, when the decision making area is split into three parts, the opposite ones (upper approximation area, and lower approximation area) and the neutral one (border area, in between the upper and lower area).

### 47. Grey System Theory as a Neutrosophication

A Grey System [30] is referring to a *grey area* (as <neutA> in neutrosophy), between extremes (as <A> and <antiA> in neutrosophy).

According to the Grey System Theory, a system with perfect information (<A>) may have a unique solution, while a system with no information (<antiA>) has no solution. In the middle (<neutA>), or grey area, of these opposite systems, there may be many available solutions (with partial information known and partial information unknown) from which an approximate solution can be extracted.

## 48. Three-Ways Decision as particular cases of Neutrosophication and of Neutrosophic Probability [31, 32, 33, 34, 35, 36]

### 48.1. Neutrosophication

Let <A> be an attribute value, <antiA> the opposite of this attribute value, and <neutA> the neutral (or indeterminate) attribute value between the opposites <A> and <antiA>.

For examples: <A> = *big*, then <antiA> = *small*, and <neutA> = *medium*; we may rewrite:

(<A>, <neutA>, <antiA>) = (big, medium, small);

or (<A>, <neutA>, <antiA>) = *(truth (*denoted as *T), indeterminacy (*denoted as *I), falsehood (*denoted as *F) )* as in Neutrosophic Logic, or (<A>, <neutA>, <antiA>) = *( membership, indeterminate-membership, monmembership )* as in Neutrosophic Set, or (<A>, <neutA>, <antiA>) = *( chance that an event occurs, indeterminate-chance that the event occurs or not, chance that the event does not occur )* as in Neutrosophic Probability, and so on.

And let's by "Concept" to mean: an item, object, idea, theory, region, universe, set, notion etc. that is characterized by this attribute.

The process of **neutrosophication** {Smarandache, 2019, [37]} means:

**a)** converting a *Classical Concept*

{ denoted as *($1_{<A>}$, $0_{<neutA>}$, $0_{<antiA>}$)-ClassicalConcept*, or *ClassicalConcept($1_{<A>}$, $0_{<neutA>}$, $0_{<antiA>}$)* }, which means that the concept is, with respect to the above attribute,

*100% <A>, 0% <neutA>, and 0% <antiA>*,

into a *Neutrosophic Concept*

{ denoted as *($T_{<A>}$, $I_{<neutA>}$, $F_{<antiA>}$)-NeutrosophicConcept*, or *NeutrosophicConcept($T_{<A>}$, $I_{<neutA>}$, $F_{<antiA>}$)* }, which means that the concept is, with respect to the above attribute,

*T% <A>, I% <neutA>, and F% <antiA>*,

which more accurately reflects our imperfect, non-idealistic reality,

where all *T, I, F* are subsets of *[0, 1]* with no other restriction;

**b)** or converting a *Fuzzy Concept,* or *Intuitionistic Fuzzy Concept* into a *Neutrosophic Concept*;

c) or converting other Concepts such as *Inconsistent Intuitionistic Fuzzy (Picture Fuzzy, Ternary Fuzzy) Concept,* or *Pythagorean Fuzzy Concept,* or *Spherical Fuzzy Concept,* or *q-Rung Orthopair Fuzzy* etc.

into a *Neutrosophic Concept* or into a *Refined Neutrosophic Concept* (i.e. $T_1\%$ <$A_1$>, $T_2\%$ <$A_2$>, ...; $I_1\%$ <$neutA_1$>, $I_2\%$ <$neutA_2$>, ...; and $F_1\%$ <$antiA_1$>, $F_2\%$ <$antiA_2$>, ...), where all $T_1, T_2, ...; I_1, I_2, ...; F_1, F_2, ...$ are subsets of *[0, 1]* with no other restriction.

d) or converting a *crisp real number* ( *r* ) into a *neutrosophic real number* of the form $r = a + bI$, where "*I*" means (literal or numerical) indeterminacy, *a* and *b* are real numbers, and "*a*" represents the determinate part of the crisp real number *r*, while *bI* the indeterminate part of *r*;

e) or converting a *crisp complex number* ( *c* ) into a *neutrosophic complex number* of the form $c = a_1 + b_1 i + (a_2 + b_2 i)I = a_1 + a_2 I + (b_1 + b_2 I)i$, where "*I*" means (literal or numerical) indeterminacy, $i = \sqrt{-1}$, with $a_1, a_2, b_1, b_2$ real numbers, and "$a_1 + b_1 i$" represents the determinate part of the complex real number *c*, while $a_2 + b_2 i$ the indeterminate part of *c*;
(we may also interpret that as: $a_1$ is the determinate part of the real-part of *c*, and $b_1$ is the determinate part of the imaginary-part of *c*; while $a_2$ is the indeterminate part of the real-part of *c*, and $b_2$ is the indeterminate part of the imaginary-part of *c*);

f) converting a *crisp, fuzzy,* or *intuitionistic fuzzy,* or *inconsistent intuitionistic fuzzy (picture fuzzy, ternary fuzzy set),* or *Pythagorean fuzzy*, or *spherical fuzzy*, or *q-rung orthopair fuzzy number* and other numbers into a *quadruple neutrosophic number* of the form $a + bT + cI + dF$, where *a, b, c, d* are real or complex numbers, while *T, I, F* are the neutrosophic components.

While the process of **deneutrosophication** means going backwards with respect to any of the above processes of neutrosophication.

Example 1:

Let the attribute <A> = cold temperature, then <antiA> = hot temperature, and <neutA> = medium temperature.

Let the concept be a country *M*, such that its northern part (30% of country's territory) is cold, its southern part is hot (50%), and in the middle there is a buffer zone with medium temperature (20%). We write:

$$M(\ 0.3_{cold\ temperature},\ 0.2_{medium\ temperature},\ 0.5_{hot\ temperature}\ )$$

where we took single-valued numbers for the neutrosophic components $T_M = 0.3$, $I_M = 0.2$, $F_M = 0.5$, and the neutrosophic components are considered dependent so their sum is equal to 1.

### 48.2. Three-Ways Decision is a particular case of Neutrosophication

Neutrosophy (based on <A>, <neutA>, <antiA>) was proposed by Smarandache [1] in 1998, and Three-Ways Decision by Yao [31] in 2009.

In Three-Ways Decision, the universe set is split into three different distinct areas, in regard to the decision process, representing:

*Acceptance, Noncommitment,* and *Rejection* respectively.

In this case, the decision attribute value $<A>$ = Acceptance, whence $<neutA>$ = Noncommitment, and $<antiA>$ = Rejection.

The classical concept = *UniverseSet*.

Therefore, we got the *NeutrosophicConcept( $T_{<A>}$, $I_{<neutA>}$, $F_{<antiA>}$ )*, denoted as:

*UniverseSet( $T_{Acceptance}$, $I_{Noncommitment}$, $F_{Rejection}$ )*,

where $T_{Acceptance}$ = *universe set's zone of acceptance*, $I_{Noncommitment}$ = *universe set's zone of noncomitment (indeterminacy)*, $F_{Rejection}$ = *= universe set's zone of rejection.*

### 48.3. Three-Ways Decision as a particular case of Neutrosophic Probability

Let's consider the event, taking a decision on a universe set.

According to Neutrosophic Probability (NP) [1, 11] one has:

*NP(decision)* = ( the universe set's elements for which the chance of the decision may be accept; the universe set's elements for which there may be an indeterminate-chance of the decision; the universe set's elements for which the chance of the decision may be reject ).

### 48.4. Refined Neutrosophy

*Refined Neutrosophy* was introduced by Smarandache [9] in 2013 and it is described as follows:

$<A>$ is refined (split) into subcomponents $<A_1>$, $<A_2>$, ..., $<A_p>$;

$<neutA>$ is refined (split) into subcomponents $<neutA_1>$, $<neutA_2>$, ..., $<neutA_r>$;

and $<antiA>$ is refined (split) into subcomponents $<antiA_1>$, $<antiA_2>$, ..., $<antiA_s>$;

where *p, r, s ≥ 1* are integers, and *p + r + s ≥ 4*.

Refined Neutrosophy is a generalization of Neutrosophy.

Example 2.

If $<A>$ = voting in country M, them $<A_1>$ = voting in Region 1 of country M for a given candidate, $<A_2>$ = voting in Region 2 of country M for a given candidate, and so on.

Similarly, $<neutA_1>$ = not voting (or casting a white or a black vote) in Region 1 of country M, $<A_2>$ = not voting in Region 2 of country M, and so on.

And $<antiA_1>$ = voting in Region 1 of country M against the given candidate, $<A_2>$ = voting in Region 2 of country M against the given candidate, and so on.

### 48.5. Extension of Three-Ways Decision to n-Ways Decision

*n-Way Decision* was introduced by Smarandache [37] in 2019.

In n-Ways Decision, the universe set is split into $n \geq 4$ different distinct areas, in regard to the decision process, representing:

*Levels of Acceptance, Levels of Noncommitment,* and *Levels of Rejection* respectively.

Levels of Acceptance may be: Very High Level of Acceptance *($<A_1>$)*, High Level of Acceptance *($<A_2>$)*, Medium Level of Acceptance *($<A_3>$)*, etc.

Similarly, Levels of Noncommitment may be: Very High Level of Noncommitment *($<neutA_1>$)*, High Level of Noncommitment *($<neutA_2>$)*, Medium Level of Noncommitment *($<neutA_3>$)*, etc.

And Levels of Rejection may be: Very High Level of Rejection *($<antiA_1>$)*, High Level of Rejection *($<antiA_2>$)*, Medium Level of Rejection *($<antiA_3>$)*, etc.

Then the *Refined Neutrosophic Concept*

{ denoted as *($T1_{<A1>}, T2_{<A2>}, …, Tp_{<Ap>}; I1_{<neutA1>}, I2_{<neutA2>}, …, Ir_{<neutAr>};$*

*$F1_{<antiA1>}, F2_{<antiA2>}, Fs_{<antiAs>}$)-RefinedNeutrosophicConcept*,

or **RefinedNeutrosophicConcept**(*$T1_{<A1>}, T2_{<A2>}, …, Tp_{<Ap>}; I1_{<neutA1>}, I2_{<neutA2>}, …, Ir_{<neutAr>}; F1_{<antiA1>}, F2_{<antiA2>}, Fs_{<antiAs>}$)*},

which means that the concept is, with respect to the above attribute value levels,

$$T1\% <A1>, T2\% <A2>, …, Tp\% <Ap>;$$

$$I1\% <neutA1>, I2\% <neutA2>, …, Ir\% <neutAr>;$$

$$F1\% <antiA1>, F2\% <antiA2>, Fs\% <antiAs>;$$

which more accurately reflects our imperfect, non-idealistic reality,

with where $p, r, s \geq 1$ are integers, and $p + r + s \geq 4$,

where all *T1, T2, ..., Tp, I1, I2, ..., Ir, F1, F2, ..., Fs* are subsets of *[0, 1]* with no other restriction.

## 49. Many More Distinctions between Neutrosophic Set (NS) and Intuitionistic Fuzzy Set (IFS) and other type sets

49.1. Neutrosophic Set can distinguish between absolute and relative

- *absolute membership* (i.e. membership in all possible worlds; we have extended Leibniz's absolute truth to absolute membership), and
- *relative membership* (membership in at least one world, but not in all), because

$$NS \text{ (absolute membership element)} = 1^+$$

while

- NS (relative membership element) = 1.

This has application in philosophy (see the *neutrosophy*). That's why the unitary standard interval $[0, 1]$ used in IFS has been extended to the unitary non-standard interval $]^-0, 1^+[$ in NS.

Similar distinctions for *absolute* or *relative non-membership*, and *absolute* or *relative indeterminate appurtenance* are allowed in NS.

While IFS cannot distinguish the absoluteness from relativeness of the components.

49.2. In NS, there is no restriction on T, I, F other than they be subsets of $]^-0, 1^+[$, thus:

$$^-0 \leq \inf T + \inf I + \inf F \leq \sup T + \sup I + \sup F \leq 3^+.$$

The inequalities (2.1) and (2.4) [17] of IFS are relaxed in NS.

This non-restriction allows paraconsistent, dialetheist, and incomplete information to be characterized in NS {i.e. the sum of all three components if they are defined as points, or sum of superior limits of all three components if they are defined as subsets can be >1 (for paraconsistent information coming from different sources), or < 1 for incomplete information}, while that information cannot be described in IFS because in IFS the components T (membership), I (indeterminacy), F (non-membership) are restricted either to $t + i + f = 1$ or to $t^2 + f^2 \leq 1$, if T, I, F are all reduced to the points (single-valued numbers) t, i, f respectively, or to $\sup T + \sup I + \sup F = 1$ if T, I, F are subsets of $[0, 1]$. Of course, there are cases when paraconsistent and incomplete informations can be normalized to 1, but this procedure is not always suitable.

In IFS paraconsistent, dialetheist, and incomplete information cannot be characterized.

This most important distinction between IFS and NS is showed in the below **Neutrosophic Cube** A'B'C'D'E'F'G'H' introduced by J. Dezert [38] in 2002.

Because in technical applications only the classical interval $[0,1]$ is used as range for the neutrosophic parameters $t, i, f$, we call the cube *ABCDEDGH* the **technical neutrosophic cube** and its extension $A'B'C'D'E'D'G'H'$ the **neutrosophic cube** (or **nonstandard neutrosophic cube**), used in the fields where we need to differentiate between *absolute* and *relative* (as in philosophy) notions.

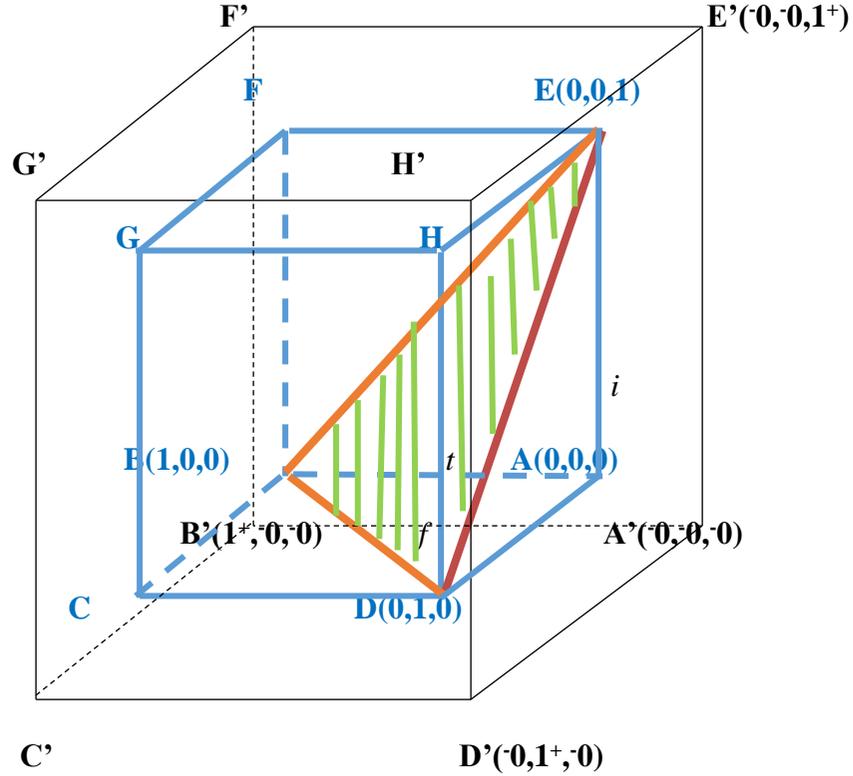

*Fig. 1. Neutrosophic Cube*

Let's consider a 3D Cartesian system of coordinates, where $t$ is the truth axis with value range in $]^-0,1^+[$, $f$ is the false axis with value range in $]^-0,1^+[$, and similarly $i$ is the indeterminate axis with value range in $]^-0,1^+[$.

We now divide the technical neutrosophic cube $ABCDEDGH$ into three disjoint regions:

a) The shaded equilateral triangle $BDE$, whose sides are equal to $\sqrt{2}$, which represents the geometrical locus of the points whose sum of the coordinates is 1.

If a point $Q$ is situated on the sides or inside of the triangle $BDE$, then $t_Q + i_Q + f_Q = 1$ as in Atanassov-intuitionistic fuzzy set $(A-IFS)$.

It is clear that IFS triangle is a restriction of (strictly included in) the NS cube.

b) The pyramid $EABD$ {situated in the right side of the $\Delta EBD$, including its faces $\Delta ABD$ (base), $\Delta EBA$, and $\Delta EDA$ (lateral faces), but excluding its face $\Delta BDE$} is the locus of the points whose sum of coordinates is less than 1.

If $P \in EABD$ then $t_P + i_P + f_P < 1$ as in inconsistent intuitionistic fuzzy set (with incomplete information).

c) In the left side of $\Delta BDE$ in the cube there is the solid $EFGCDEBD$ (excluding $\Delta BDE$) which is the locus of points whose sum of their coordinates is greater than 1 as in the paraconsistent set.

If a point $R \in EFGCDEBD$, then $t_R + i_R + f_R > 1$.

It is possible to get the **sum of coordinates strictly less than 1 or strictly greater than 1**. For example having three independent sources of information:

- We have a source which is capable to find only the degree of membership of an element; but it is unable to find the degree of non-membership;
- Another source which is capable to find only the degree of non-membership of an element;
- Or a source which only computes the indeterminacy.

Thus, when we put the results together of these sources, it is possible that their sum is not 1, but smaller or greater.

Also, in information fusion, when dealing with indeterminate models (i.e. elements of the fusion space which are indeterminate/unknown, such as intersections we don't know if they are empty or not since we don't have enough information, similarly for complements of indeterminate elements, etc.): if we compute the believe in that element (truth), the disbelieve in that element (falsehood), and the indeterminacy part of that element, then the sum of these three components is strictly less than 1 (the difference to 1 is the missing information).

49.3) Relation (2.3) from interval-valued intuitionistic fuzzy set is relaxed in NS, i.e. the intervals do not necessarily belong to Int[0,1] but to [0,1], even more general to $]^-0, 1^+[$.

49.4) In NS the components T, I, F can also be ***nonstandard subsets*** included in the unitary non-standard interval $]^-0, 1^+[$, not only *standard* subsets included in the unitary standard interval [0, 1] as in IFS.

49.5) NS, like dialetheism, can describe **paradoxist elements**, NS(paradoxist element) = (1, 1, 1), while IFL cannot describe a paradox because the sum of components should be 1 in IFS.

49.6) The connectors/operators in IFS are defined with respect to T and F only, i.e. membership and nonmembership only (hence the Indeterminacy is what's left from 1), while in NS they can be defined with respect to any of them (no restriction).

But, for interval-valued intuitionistic fuzzy set one cannot find any left indeterminacy.

49.7) Component "*I*", indeterminacy, can be split into more subcomponents in order to better catch the vague information we work with, and such, for example, one can get more accurate answers to the *Question-Answering Systems* initiated by Zadeh (2003).

{In Belnap's four-valued logic (1977) indeterminacy is split into Uncertainty (*U*) and Contradiction (*C*), but they were interrelated.}

Even more, one can split "I" into Contradiction, Uncertainty, and Unknown, and we get a five-valued logic.

In a general Refined Neutrosophic Logic, T can be split into subcomponents $T_1, T_2, ..., T_p$, and I into $I_1, I_2, ..., I_r$, and F into $F_1, F_2, ...,F_s$, where $p, r, s \geq 1$ and $p + r + s = n \geq 3$. Even more: T, I, and/or F (or any of their subcomponents $T_j, I_k$, and/or $F_l$) can be countable or uncountable infinite sets.

49.8) Indeterminacy is independent from membership/truth and non-membership/falsehood in NS/Nl, while in IFS/IFL it is not.
In neutrosophics there are two types of indeterminacies:

> a) *Numerical Indeterminacy (or Degree of Indeterminacy)*, which has the form $(t, i, f) \neq (1, 0, 0)$, where *t, i, f* are numbers, intervals, or subsets included in the unit interval [0, 1], and it is the base for the *(t, i, f)-Neutrosophic Structures*.

> b) *Non-numerical Indeterminacy (or Literal Indeterminacy)*, which is the letter "*I*" standing for unknown (non-determinate), such that $I^2 = I$, and used in the composition of the neutrosophic number $N = a + bI$, where *a* and *b* are real or complex numbers, and *a* is the determinate part of number *N*, while *bI* is the indeterminate part of *N*. The neutrosophic numbers are the base for the *I*-Neutrosophic Structures.

49.9) NS has a better and clear terminology (name) as "neutrosophic" (which means the neutral part: i.e. neither true/membership nor false/nonmembership), while IFS's name "intuitionistic" produces confusion with Intuitionistic Logic, which is something different (see the article by Didier Dubois et al. [39], 2005).

49.10) The Neutrosophic Set was extended [Smarandache, 2007] to **Neutrosophic Overset** (when some neutrosophic component is > 1), and to **Neutrosophic Underset** (when some neutrosophic component is < 0), and to and to **Neutrosophic Offset** (when some neutrosophic components are off the interval [0, 1], i.e. some neutrosophic component > 1 and some neutrosophic component < 0). In IFS the degree of a component is not allowed to be outside of the classical interval [0, 1].

This is no surprise with respect to the classical fuzzy set/logic, intuitionistic fuzzy set/logic, or classical and imprecise probability where the values are not allowed outside the interval [0, 1], since our real-world has numerous examples and applications of over/under/off neutrosophic components.

Example:

In a given company a full-time employer works 40 hours per week. Let's consider the last week period.

Helen worked part-time, only 30 hours, and the other 10 hours she was absent without payment; hence, her membership degree was $30/40 = 0.75 < 1$.

John worked full-time, 40 hours, so he had the membership degree $40/40 = 1$, with respect to this company.

But George worked overtime 5 hours, so his membership degree was (40+5)/40 = 45/40 = 1.125 > 1. Thus, we need to make distinction between employees who work overtime, and those who work full-time or part-time. That's why we need to associate a degree of membership greater than 1 to the overtime workers.

Now, another employee, Jane, was absent without pay for the whole week, so her degree of membership was 0/40 = 0.

Yet, Richard, who was also hired as a full-time, not only didn't come to work last week at all (0 worked hours), but he produced, by accidentally starting a devastating fire, much damage to the company, which was estimated at a value half of his salary (i.e. as he would have gotten for working 20 hours). Therefore, his membership degree has to be less that Jane's (since Jane produced no damage). Whence, Richard's degree of membership with respect to this company was - 20/40 = - 0.50 < 0.

Therefore, the membership degrees > 1 and < 0 are real in our world, so we have to take them into consideration.

Then, similarly, the Neutrosophic Logic/Measure/Probability/Statistics etc. were extended to respectively Neutrosophic Over/Under/Off Logic, Measure, Probability, Statistics etc. {Smarandache, 2007 [8]}.

49.11) **Neutrosophic Tripolar** (and in general **Multipolar**) **Set** and **Logic** {Smarandache, 2007 [8]} of the form:

( <$T^+_1, T^+_2, …, T^+_n; T^0; T^-_{-n}, …, T^-_{-2}, T^-_{-1}$>, <$I^+_1, I^+_2, …, I^+_n; I^0; I^-_{-n}, …, I^-_{-2}, I^-_{-1}$>,

<$F^+_1, F^+_2, …, F^+_n; F^0; F^-_{-n}, …, F^-_{-2}, F^-_{-1}$> )

where we have multiple positive/neutral/negative degrees of T, I, and F respectively.

49.12) The **Neutrosophic Numbers** have been introduced by W.B. Vasantha Kandasamy and F. Smarandache [40] in 2003, which are numbers of the form N = *a* + *bI*, where *a, b* are real or complex numbers, while "I" is the indeterminacy part of the neutrosophic number N, such that $I^2$ = I and $αI+βI = (α+β)I$.

Of course, indeterminacy "I" is different from the imaginary unit $i = \sqrt{-1}$.
In general one has $I^n$ = I if n > 0, and it is undefined if n ≤ 0.
49.13) Also, **Neutrosophic Refined Numbers** were introduced (Smarandache [31], 2015) as:
*a* + *b₁I₁* + *b₂I₂* + ... + *bₘIₘ*, where *a, b₁, b₂, …, bₘ* are real or complex numbers, while the *I₁, I₂, …, Iₘ* are types of sub-indeterminacies, for *m ≥ 1*.
49.14) The algebraic structures using neutrosophic numbers gave birth to the ***I*-Neutrosophic Algebraic Structures** [see for example "neutrosophic groups", "neutrosophic rings", "neutrosophic vector space", "neutrosophic matrices, bimatrices, …, n-matrices", etc.], introduced by W.B. Vasantha Kandasamy, Ilanthenral K., F. Smarandache [41] et al. since 2003.

Example of Neutrosophic Matrix: $\begin{bmatrix} 1 & 2+I & -5 \\ 0 & 1/3 & I \\ -1+4I & 6 & 5I \end{bmatrix}$.

Example of Neutrosophic Ring: ({a+bI, with a, b ϵ R}, +, ·), where of course (a+bI)+(c+dI) = (a+c)+(b+d)I, and (a+bI) · (c+dI) = (ac) + (ad+bc+bd)I.

49.15) Also, to **Refined *I*-Neutrosophic Algebraic Structures**, which are structures using sets of refined neutrosophic numbers [41].

49.16) Types of Neutrosophic Graphs (and Trees):
a-c) Indeterminacy "I" led to the definition of the **Neutrosophic Graphs** (graphs which have: either at least one indeterminate edge, or at least one indeterminate vertex, or both some indeterminate edge and some indeterminate vertex), and **Neutrosophic Trees** (trees which have: either at least one indeterminate edge, or at least one indeterminate vertex, or both some indeterminate edge and some indeterminate vertex), which have many applications in social sciences.
Another type of neutrosophic graph is when at least one edge has a neutrosophic (t, i, f) truth-value. As a consequence, the Neutrosophic Cognitive Maps (Vasantha & Smarandache, 2003]) and Neutrosophic Relational Maps (Vasantha & Smarandache, 2004) are generalizations of fuzzy cognitive maps and respectively fuzzy relational maps, Neutrosophic Relational Equations (Vasantha & Smarandache, 2004), Neutrosophic Relational Data (Wang, Smarandache, Sunderraman, Rogatko - 2008), etc.
A Neutrosophic Cognitive Map is a neutrosophic directed graph with concepts like policies, events etc. as vertices, and causalities or indeterminates as edges. It represents the causal relationship between concepts.

An edge is said indeterminate if we don't know if it is any relationship between the vertices it connects, or for a directed graph we don't know if it is a directly or inversely proportional relationship. We may write for such edge that (t, i, f) = (0, 1, 0).
A vertex is indeterminate if we don't know what kind of vertex it is since we have incomplete information. We may write for such vertex that (t, i, f) = (0, 1, 0).

Example of Neutrosophic Graph (edges $V_1V_3$, $V_1V_5$, $V_2V_3$ are indeterminate and they are drawn as dotted):

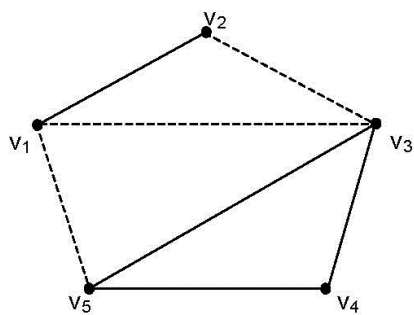

*Fig. 2. Neutrosophic Graph { with I (indeterminate) edges }*

and its neutrosophic adjacency matrix is:

$$\begin{bmatrix} 0 & 1 & I & 0 & I \\ 1 & 0 & I & 0 & 0 \\ I & I & 0 & 1 & 1 \\ 0 & 0 & 1 & 0 & 1 \\ I & 0 & 1 & 1 & 0 \end{bmatrix}$$

*Fig. 3. Neutrosophic Adjacency Matrix of the Neutrosophic Graph*

The edges mean: 0 = no connection between vertices, 1 = connection between vertices, I = indeterminate connection (not known if it is, or if it is not).

Such notions are not used in the fuzzy theory.

Example of Neutrosophic Cognitive Map (NCM), which is a generalization of the Fuzzy Cognitive Maps.

Let's have the following vertices:
$C_1$ - Child Labor
$C_2$ - Political Leaders
$C_3$ - Good Teachers
$C_4$ - Poverty
$C_5$ - Industrialists
$C_6$ - Public practicing/encouraging Child Labor
$C_7$ - Good Non-Governmental Organizations (NGOs)

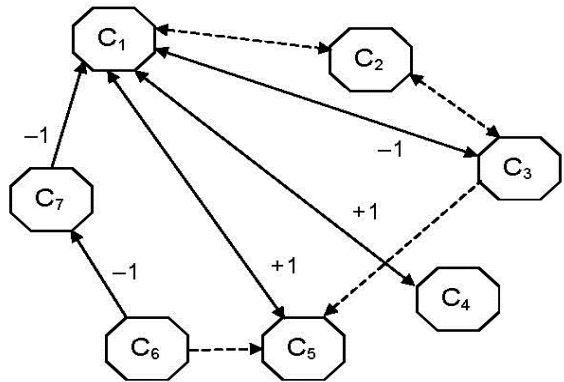

*Fig. 4. Neutrosophic Cognitive Map*

The corresponding neutrosophic adjacency matrix related to this neutrosophic cognitive map is:

$$\begin{bmatrix} 0 & I & -1 & 1 & 1 & 0 & 0 \\ I & 0 & I & 0 & 0 & 0 & 0 \\ -1 & I & 0 & 0 & I & 0 & 0 \\ 1 & 0 & 0 & 0 & 0 & 0 & 0 \\ 1 & 0 & 0 & 0 & 0 & 0 & 0 \\ 0 & 0 & 0 & 0 & I & 0 & -1 \\ -1 & 0 & 0 & 0 & 0 & 0 & 0 \end{bmatrix}$$

*Fig. 4. Neutrosophic Adjacency Matrix of the Neutrosophic Cognitive Map*

The edges mean: 0 = no connection between vertices, 1 = directly proportional connection, -1 = inversely proportionally connection, and I = indeterminate connection (not knowing what kind of relationship is between the vertices that the edge connects).
Such **literal indeterminacy** (letter **I**) does not occur in previous set theories, including intuitionistic fuzzy set; they had only *numerical indeterminacy*.

d) Another type of neutrosophic graphs (and trees) [Smarandache, 2015, [41]]:
An edge of a graph, let's say from A to B (i.e. how A influences B),
may have a neutrosophic value *(t, i, f)*,
where t means the positive influence of A on B,
    i means the indeterminate influence of A on B,
  and f means the negative influence of A on B.

Then, if we have, let's say: A->B->C such that A->B has the neutrosophic value $(t_1, i_1, f_1)$ and B->C has the neutrosophic value $(t_2, i_2, f_2)$, then A->C has the neutrosophic value $(t_1, i_1, f_1) \wedge (t_2, i_2, f_2)$, where $\wedge$ is the AND neutrosophic operator.

e) Also, again a different type of graph: we can consider a vertex A as: t% belonging/membership to the graph, i% indeterminate membership to the graph, and f% nonmembership
to the graph.
13.f) Any of the previous types of graphs (or trees) put together.
13.g) **Tripolar** (and **Multipolar**) **Graph**, which is a graph whose vertexes or edges have the form $(<T^+, T^0, T^->, <I^+, I^0, I^->, <F^+, F^0, F^->)$ and respectively: $(<T^+_j, T^0, T^-_j>, <I^+_j, I^0, I^-_j>, <F^+_j, F^0, F^-_j>)$.

49.17) The **Neutrosophic Probability** (NP), introduced in 1995, was extended and developed as a generalization of the classical and imprecise probabilities {Smarandache, 2013 [11]}. NP of an event $\mathscr{E}$ is the chance that event $\mathscr{E}$ occurs, the chance that event $\mathscr{E}$ doesn't occur, and the chance of indeterminacy (not knowing if the event $\mathscr{E}$ occurs or not).

In classical probability $n_{sup} \leq 1$, while in neutrosophic probability $n_{sup} \leq 3^+$.

In imprecise probability: the probability of an event is a subset T in [0, 1], not a number p in [0, 1], what's left is supposed to be the opposite, subset F (also from the unit interval [0, 1]); there is no indeterminate subset I in imprecise probability.

In neutrosophic probability one has, besides randomness, indeterminacy due to construction materials and shapes of the probability elements and space.
In consequence, neutrosophic probability deals with two types of variables: random variables and indeterminacy variables, and two types of processes: stochastic process and respectively indeterminate process.

49.18) And consequently the **Neutrosophic Statistics**, introduced in 1995 and developed in {Smarandache, 2014, [12]}, which is the analysis of the neutrosophic events.
Neutrosophic Statistics means statistical analysis of population or sample that has indeterminate (imprecise, ambiguous, vague, incomplete, unknown) data. For example, the population or sample size might not be exactly determinate because of some individuals that partially belong to the population or sample, and partially they do not belong, or individuals whose appurtenance is completely unknown. Also, there are population or sample individuals whose data could be indeterminate. It is possible to define the neutrosophic statistics in many ways, because there are various types of indeterminacies, depending on the problem to solve.

Neutrosophic statistics deals with neutrosophic numbers, neutrosophic probability distribution, neutrosophic estimation, neutrosophic regression.

The function that models the neutrosophic probability of a random variable x is called *neutrosophic distribution*: $NP(x) = (T(x), I(x), F(x))$, where $T(x)$ represents the probability that value x occurs, $F(x)$ represents the probability that value x does not occur, and $I(x)$ represents the indeterminate / unknown probability of value x.

49.19) Also, **Neutrosophic Measure** and **Neutrosophic Integral** were introduced {Smarandache, 2013, [11]}.

49.20) **Neutrosophy** {Smarandache, 1995, [1, 2, 3, 4, 5, 7]} opened a new field in philosophy.

<u>Neutrosophy</u> is a new branch of philosophy that studies the origin, nature, and scope of neutralities, as well as their interactions with different ideational spectra.

This theory considers every notion or idea <A> together with its opposite or negation <Anti-A> and the spectrum of "neutralities" <Neut-A> (i.e. notions or ideas located between the two extremes, supporting neither <A> nor <Anti-A>). The <Neut-A> and <Anti-A> ideas together are referred to as <Non-A>.

According to this theory every idea <A> tends to be neutralized and balanced by <Anti-A> and <Non-A> ideas - as a state of equilibrium.

In a classical way <A>, <Neut-A>, <Anti-A> are disjoint two by two.

But, since in many cases the borders between notions are vague, imprecise, Sorites, it is possible that <A>, <Neut-A>, <Anti-A> (and <Non-A> of course) have common parts two by two as well.

Neutrosophy is the base of neutrosophic logic, neutrosophic set, neutrosophic probability and statistics used in engineering applications (especially for software and information fusion), medicine, military, cybernetics, physics.

We have extended dialectics (based on the opposites <A> and <antiA>) to neutrosophy (based on <A>, <antiA> and <neutA>.

49.21) In consequence, we extended the thesis-antithesis-synthesis to thesis-antithesis-neutrothesis-neutrosynthesis {Smarandache, 2015 [41]}.

49.22) Neutrosophy extended the Law of Included Middle to the **Law of Included Multiple-Middle** {Smarandache, 2014 [10]} in accordance with the n-valued refined neutrosophic logic.

49.23) Smarandache (2015 [41]) introduced the **Neutrosophic Axiomatic System** and **Neutrosophic Deducibility**.

49.24) Then he introduced the **(t, i, f)-Neutrosophic Structure** (2015 [41]), which is a structure whose space, or at least one of its axioms (laws), has some indeterminacy of the form $(t, i, f) \neq (1, 0, 0)$.

Also, we defined the combined *(t, i, f)-I-Neutrosophic Algebraic Structures*, i.e. algebraic structures based on neutrosophic numbers of the form $a + bI$, but also having some indeterminacy [ of the form $(t, i, f) \neq (1, 0, 0)$ ] related to the structure space (i.e. elements which only partially belong to the space, or elements we know nothing if they belong to the space or not) or indeterminacy [ of the form $(t, i, f) \neq (1, 0, 0)$ ] related to at least one axiom (or law) acting on the structure space) .
    Even more, we generalized them to *Refined (t, i, f)- Refined I-Neutrosophic Algebraic Structures,* or *$(t_j, i_k, f_l)$-$i_s$-Neutrosophic Algebraic Structures*; where $t_j$ means that $t$ has been refined to $j$ subcomponents $t_1, t_2, …, t_j$; similarly for $i_k, f_l$ and respectively $i_s$.

49.25) Smarandache and Ali [2014-2016] introduced the *Neutrosophic Triplet Structures* [42, 43, 44].

A *Neutrosophic Triplet*, is a triplet of the form:
$$< a, neut(a), anti(a) >,$$
where *neut(a)* is the neutral of a, i.e. an element (different from the identity element of the operation *) such that *a*neut(a) = neut(a)*a = a,*
while anti(a) is the opposite of a, i.e. an element such that *a*anti(a) = anti(a)*a = neut(a).*
Neutrosophy means not only indeterminacy, but also neutral (i.e. neither true nor false).
For example we can have neutrosophic triplet semigroups, neutrosophic triplet loops, etc.

Further on Smaradnache extended the neutrosophic triplet < *a, neut(a), anti(a)* > to a
                   **$m$-valued refined neutrosophic triplet**,

in a similar way as it was done for $T_1, T_2, ...; I_1, I_2, ...; F_1, F_2, ...$ (i.e. the refinement of neutrosophic components).

It will work in some cases, depending on the composition law *. It depends on each * how many neutrals and anti's there is for each element "a".

We may have an m-tuple with respect to the element "a" in the following way:
( $a; neut_1(a), neut_2(a), ..., neut_p(a); anti_1(a), anti_2(a), ..., anti_p(a)$ ),
where $m = 1+2p$,
such that:
- all $neut_1(a), neut_2(a), ..., neut_p(a)$ are distinct two by two, and each one is different from the unitary element with respect to the composition law *;
- also:

$$a*neut_1(a) = neut_1(a)*a = a$$
$$a*neut_2(a) = neut_2(a)*a = a$$
$$...........................$$
$$a*neut_p(a) = neut_p(a)*a = a;$$
- and
$$a*anti_1(a) = anti_1(a)*a = neut_1(a)$$
$$a*anti_2(a) = anti_2(a)*a = neut_2(a)$$
$$...............................$$
$$a*anti_p(a) = anti_p(a)*a = neut_p(a);$$

- where all $anti_1(a), anti_2(a), ..., anti_p(a)$ are distinct two by two, and in case when there are duplicates, the duplicates are discarded.

49.26) As latest minute development, the crisp, fuzzy, intuitionistic fuzzy, inconsistent intuitionistic fuzzy (picture fuzzy, ternary fuzzy), and neutrosophic sets were extended by Smarandache [45] in 2017 to **plithogenic set**, which is:

A set *P* whose elements are characterized by many attributes' values. An attribute value v has a corresponding (fuzzy, intuitionistic fuzzy, picture fuzzy, or neutrosophic) degree of appurtenance *d(x,v)* of the element *x*, to the set *P*, with respect to some given criteria. In order to obtain a better accuracy for the *plithogenic aggregation operators* in the plithogenic set, and for a more exact inclusion (partial order), a (fuzzy, intuitionistic fuzzy, picture fuzzy, or neutrosophic) *contradiction (dissimilarity) degree* is defined between each attribute value and the dominant (most important) attribute value. The plithogenic intersection and union are *linear combinations of the fuzzy operators t-norm and t-conorm*, while the plithogenic complement (negation), inclusion (inequality), equality (equivalence) are influenced by the attribute values contradiction (dissimilarity) degrees.

### 35. Conclusion

In this paper we proved that neutrosophic set is a generalization of intuitionistic fuzzy set and inconsistent intuitionistic fuzzy set (picture fuzzy set, ternary fuzzy set).

By transforming (restraining) the neutrosophic components into intuitionistic fuzzy components, as Atanassov and Vassiliev proposed, the *independence of the components is lost* and the *indeterminacy is ignored by the intuitionistic fuzzy aggregation operators*. Also, the result after applying the neutrosophic operators is different from the result obtained after applying the intuitionistic fuzzy operators (with respect to the same problem to solve).

We presented many distinctions between neutrosophic set and intuitionistic fuzzy set, and we showed that neutrosophic set is more general and more flexible than previous set theories. Neutrosophy's applications in various fields such neutrosophic probability, neutrosophic statistics, neutrosophic algebraic structures and so on were also listed {see also [46]}.

Neutrosophic Set (NS) is also a generalization of Inconsistent Intuitionistic Fuzzy Set (IIFS) { which is equivalent to the Picture Fuzzy Set (PFS) and Ternary Fuzzy Set (TFS) }, Pythagorean Fuzzy Set (PyFS) {Atanassov's Intuitionistic Fuzzy Set of second type}, Spherical Fuzzy Set (SFS), n-HyperSpherical Fuzzy Set (n-HSFS), and q-Rung Orthopair Fuzzy Set (q-ROFS). And Refined Neutrosophic Set (RNS) is an extension of Neutrosophic Set. And all these sets are more general than Intuitionistic Fuzzy Set.

Neutrosophy is a particular case of Refined Neutrosophy, and consequently Neutrosophication is a particular case of Refined Neutrosophication. Also, Regret Theory, Grey System Theory, and Three-Ways Decision are particular cases of Neutrosophication and of Neutrosophic Probability. We have extended the Three-Ways Decision to n-Ways Decision, which is a particular case of Refined Neutrosophy.

## Acknowledgement

The author deeply thanks Dr. Said Broumi for revealing Atanassov and Vassiliev's paper [6] and for his comments.## References

[1] Florentin Smarandache, *Definition of Neutrosophic Logic – A Generalization of the Intuitionistic Fuzzy Logic*, Proceedings of the Third Conference of the European Society for Fuzzy Logic and Technology, EUSFLAT 2003, September 10-12, 2003, Zittau, Germany; University of Applied Sciences at Zittau/Goerlitz, 141-146.

[2-3-4] Florentin Smarandache, *Neutrosophic Set, A Generalization of the Intuitionistic Fuzzy Set*
  a. in International Journal of Pure and Applied Mathematics, Vol. 24, No. 3, 287-297, 2005;
  b. also in Proceedings of 2006 IEEE International Conference on Granular Computing, edited by Yan-Qing Zhang and Tsau Young Lin, Georgia State University, Atlanta, pp. 38-42, 2006;
  c. second version in Journal of Defense Resources Management, Brasov, Romania, No. 1, 107-116, 2010.

[5] Florentin Smarandache*, A Geometric Interpretation of the Neutrosophic Set – A Generalization of the Intuitionistic Fuzzy Set*, 2011 IEEE International Conference on Granular Computing, edited by Tzung-Pei Hong, Yasuo Kudo, Mineichi Kudo, Tsau-Young Lin, Been-Chian Chien, Shyue-Liang Wang, Masahiro Inuiguchi, GuiLong Liu, IEEE Computer Society, National University of Kaohsiung, Taiwan, 602-606, 8-10 November 2011, http://fs.unm.edu/IFS-generalized.pdf46

[6] Krassimir Atanassov and Peter Vassilev, *Intuitionistic fuzzy sets and other fuzzy sets extensions representable by them*, Journal of Intelligent & Fuzzy Systems, xx (20xx) x–xx, 2019, IOS Press, DOI: 10.3233/JIFS-179426 (under press).

[7] Florentin Smarandache, *Neutrosophy, A New Branch of Philosophy*, Multiple-Valued Logic / An International Journal, Vol. 8, No. 3, 297-384, 2002. This whole issue of this journal is dedicated to Neutrosophy and Neutrosophic Logic, http://fs.unm.edu/Neutrosophy-A-New-Branch-of-Philosophy.pdf

[8] Florentin Smarandache, *Neutrosophic Overset, Neutrosophic Underset, and Neutrosophic Offset. Similarly for Neutrosophic Over-/Under-/Off- Logic, Probability, and Statistics*, 168 p., Pons Editions, Bruxelles, Belgique, 2016; https://arxiv.org/ftp/arxiv/papers/1607/1607.00249.pdf

[9] Florentin Smarandache, *n-Valued Refined Neutrosophic Logic and Its Applications in Physics*, Progress in Physics, 143-146, Vol. 4, 2013; https://arxiv.org/ftp/arxiv/papers/1407/1407.1041.pdf

[10] Florentin Smarandache, *Law of Included Multiple-Middle & Principle of Dynamic Neutrosophic Opposition*, by Florentin Smarandache, EuropaNova & Educational, Brussels-Columbus (Belgium-USA), 136 p., 2014; http://fs.unm.edu/LawIncludedMultiple-Middle.pdf

[11] Florentin Smarandache, *Introduction to Neutrosophic Measure, Neutrosophic Integral and Neutrosophic Probability*, Sitech, 2013, http://fs.gallup.unm.edu/NeutrosophicMeasureIntegralProbability.pdf

[12] Florentin Smarandache, *Introduction to Neutrosophic Statistics*, Sitech, 2014, http://fs.gallup.unm.edu/NeutrosophicStatistics.pdf.

[13] Florentin Smarandache, *Neutrosophy. Neutrosophic Probability, Set, and Logic*. ProQuest Information and Learning, Ann Arbor, Michigan, USA, 105 p., 1998, 2000, 2002, 2005, 2006.

[14] Florentin Smarandache, *A Unifying Field in Logics: Neutrosophic Logic*, Multiple-Valued Logic / An International Journal, Vol. 8, No. 3, 385-438, 2002, http://fs.unm.edu /eBook-Neutrosophics6.pdf

[15] Florentin Smarandache, editor, *Proceedings of the First International Conference on Neutrosophy, Neutrosophic Logic, Neutrosophic Set, Neutrosophic Probability and Statistics*, University of New Mexico, Gallup Campus, Xiquan, Phoenix, 147 p., 2002, http://fs.unm.edu /NeutrosophicProceedings.pdf.

[16] Krassimir Atanassov, *Intuitionistic fuzzy sets*, Fuzzy Sets and Systems, 20: 87-96, 1986.

[17] Krassimir T. Atanassov, *Intuitionistic Fuzzy Sets*, Physica-Verlag, Heidelberg, N.Y., 1999.

[18] Krassimir Atanassov, *Intuitionistic fuzzy sets*, VII ITKR's Session, Sofia, June 1983 (Deposed in Central Sci. - Techn. Library of Bulg. Acad. of Sci., 1697/84) (in Bulgarian); Reprinted: Int J Bioautomation **20**(S1) (2016), pp. S1-S6.

[19] C. Hinde and R. Patching, *Inconsistent intuitionistic fuzzy sets. Developments in Fuzzy Sets, Intuitionistic Fuzzy Sets*, Generalized Nets and Related Topics 1 (2008), 133–153.

[20] B.C. Cuong and V. Kreinovich, *Picture fuzzy sets - a new concept for computational intelligence problems*, Proceedings of the Third World Congress on Information and Communication Technologies WICT'2013, Hanoi, Vietnam, December *15-18*, 2013, pp. 1-6.

[21] Chao Wang, Minghu Ha and Xiaqowei Liu, *A mathematical model of ternary fuzzy set for voting*, Journal of Intelligent & Fuzzy Systems, 29 (2015), 2381-2386.

[22] Florentin Smarandache, *Degree of Dependence and Independence of the (Sub)Components of Fuzzy Set and Neutrosophic Set*, Neutrosophic Sets and Systems, vol. 11, 2016, pp. 95-97, doi.org/10.5281/zenodo.571359,
http://fs.unm.edu/NSS/DegreeOfDependenceAndIndependence.pdf

[23] Princy R, Mohana K, *Spherical Bipolar Fuzzy Sets and its Application in Multi Criteria Decision Making Problem*, Journal of New Theory, 2019 (under press).

[24] R. R. Yager, *Pythagorean fuzzy subsets*, Joint IFSA World Congress and NAFIPS Annual Meeting (IFSA/NAFIPS), Edmonton, AB, Canada, pp. 57-61.

[25] Fatma Kutlu Gündoğdu, Cengiz Kahraman, *A novel spherical fuzzy QFD method and its application to the linear delta robot technology development*, Engineering Applications of Artificial Intelligence, 87 (2020) 103348.

[26] Abhishek Guleria, Rakesh Kumar Bajaj, *T-Spherical Fuzzy Graphs: Operations and Applications in Various Selection Processes*, Arabian Journal for Science and Engineering, https://doi.org/10.1007/s13369-019-04107-y

[27] Florentin Smarandache, *Neutrosophic Perspectives: Triplets, Duplets, Multisets, Hybrid Operators, Modal Logic, Hedge Algebras. And Applications*, Second extended and improved edition, Pons Publishing House Brussels, 2017, http://fs.unm.edu/NeutrosophicPerspectives-ed2.pdf

[28] Ronald R. Yager, *Generalized orthopair fuzzy sets*, IEEE Trans. Fuzzy Syst., vol. 25, no. 5, pp. 1222-1230, Oct. 2017.

[29] H. Bleichrodt, A. Cillo, E. Diecidue, *A quantitative measurement of regret theory*, Manage. Sci. 56(1) (2010) 161-175.

[30] J. L. Deng, *Introduction to Grey System Theory*, The Journal of Grey System, 1(1): 1-24, 1989.

[31] Y. Yao, *Three-way decision: an interpretation of rules in rough set theory*, in Proceeding of 4th International Conference on Rough Sets and Knowledge Technology, LNAI, Vol. 5589, Springer Berlin Heidelberg, 2009, pp. 642–649.

[32] Florentin Smarandache, *Three-Ways Decision is a particular case of Neutrosophication*, in volume *Nidus Idearum. Scilogs, VII: superluminal physics*, *vol. vii*, Pons Ed., Brussels, pp. 97 - 102, 2019.

[33] Prem Kumar Singh, *Three-way fuzzy concept lattice representation using neutrosophic set*, International Journal of Machine Learning and Cybernetics, 2017, Vol 8, Issue 1, pp. 69-79.

[34] Prem Kumar Singh, *Three-way n-valued neutrosophic concept lattice at different granulation*, International Journal of Machine Learning and Cybernetics, November 2018, Vol 9, Issue 11, pp. 1839-1855.

[35] Prem Kumar Singh, *Complex neutrosophic concept lattice and its applications to Air quality analysis*, Chaos, Solitons and Fractals, Elsevier, 2018, Vol 109, pp. 206-213.

[36] Prem Kumar Singh, *Interval-valued neutrosophic graph representation of concept lattice and its (α, β, γ)-decomposition*, Arabian Journal for Science and Engineering, Springer, 2018, Vol. 43, Issue 2, pp. 723-740.

[37] Florentin Smarandache, *Extension of Three-Ways Decision to n-Ways Decision*, in NIDUS IDEARUM. scilogs, VII: superluminal physics Brussels, 2019, http://fs.unm.edu/NidusIdearum7.pdf

[38] Jean Dezert, *Open Questions to Neutrosophic Inferences,* Multiple-Valued Logic / An International Journal, Vol. 8, No. 3, 439-472, June 2002.

[39] Didier Dubois, S. Gottwald, P. Hajek, Henry Prade, *Terminological difficulties in fuzzy set theory – The case of intuitionistic fuzzy sets*, Fuzzy Sets and Systems, 156 (2005), 585-491.

[40] W. B. Vasantha Kandasamy, Florentin Smarandache, *Fuzzy Cognitive Maps and Neutrosophic Cognitive Maps*, ProQuest Information & Learning, Ann Arbor, Michigan, USA, 2003.

[41] Florentin Smarandache, *Symbolic Neutrosophic Theory*, Europa Nova, Bruxelles, 194 p., 2015; https://arxiv.org/ftp/arxiv/papers/1512/1512.00047.pdf

[42] Florentin Smarandache, Mumtaz Ali, *The Neutrosophic Triplet Group and its Application to Physics*, presented by F. S. to Universidad Nacional de Quilmes, Department of Science and Technology, Bernal, Buenos Aires, Argentina, 02 June 2014.

[43] Florentin Smarandache and Mumtaz Ali, *Neutrosophic Triplet Group*, Neural Computing and Applications, Springer, 1-7, 2016, DOI: 10.1007/s00521-016-2535-x, http://fs.unm.edu/NeutrosophicTriplets.htm

[44] F. Smarandache, M. Ali, *Neutrosophic triplet as extension of matter plasma, unmatter plasma, and antimatter plasma*, 69th annual gaseous electronics conference, Bochum, Germany, Veranstaltungszentrum & Audimax, Ruhr-Universitat, 10–14 Oct. 2016, http://meetings.aps.org/Meeting/GEC16/Session/HT6.111

[45] Florentin Smarandache, *Plithogeny, Plithogenic Set, Logic, Probability, and Statistics*, Infinite Study Publ. Hse., GoogleLLC, Mountain View, California, USA, 2017, https://arxiv.org/ftp/arxiv/papers/1808/1808.03948.pdf Harvard SAO/NASA ADS: http://adsabs.harvard.edu/cgi-bin/bib_query?arXiv:1808.03948

[46] Florentin Smarandache, *Neutrosophic Set as Generalization of Intuitioonistic Fuzzy Set, Picture Fuzzy Set and Spherical Fuzzy Set, and its Physical Applications*, 2019 Joint Fall Meeting of the Texas Sections of American Physical Society (APS), AAPT and Zone 13 of the SPS, Friday–Saturday, October 25–26, 2019; Lubbock, Texas, USA, http://meetings.aps.org/Meeting/TSF19/scheduling?ukey=1480464-TSF19-otUQvu

*This version is slightly updated with respect to the previously published versions. ]*